\documentclass[11pt, reqno]{amsart}
\usepackage{etex}
\usepackage{textcmds} 
\usepackage{amsmath, amssymb, amsfonts, amstext, verbatim, amsthm, mathrsfs}
\usepackage{microtype}
\usepackage[all]{xy}
\usepackage[usenames,dvipsnames]{xcolor}
\usepackage{aliascnt}
\usepackage{enumitem}
\usepackage{xspace}
\usepackage{rotating}
\usepackage[margin=3.5cm]{geometry}
\usepackage{dsfont}
\usepackage{bm}
\usepackage{color}
\usepackage{subfigure}
\usepackage{array}
\usepackage{euscript}
\usepackage[T1]{fontenc}
\usepackage{mathbbol}
\usepackage{nccmath}
\usepackage{marginnote}
\usepackage{stmaryrd}
\usepackage{youngtab}
\usepackage{tikz,pgfplots}
\usetikzlibrary{calc}
\usepackage{blkarray}
\usepackage[normalem]{ulem}
\usepackage[abs]{overpic}
\usepackage{cite}
\usepackage[pagebackref,colorlinks=true,linkcolor=black,citecolor=blue,urlcolor=blue,citebordercolor={0 0 1},urlbordercolor={0 0 1},linkbordercolor={0 0 1}]{hyperref}

\renewcommand*{\backrefalt}[4]{%
  \ifcase #1 
  \or (Cited on page~#2.)%
  \else (Cited on pages~#2.)%
  \fi}

\usepackage{amsmath}
\usetikzlibrary{arrows.meta}
\definecolor{volume}{HTML}{5260D8}
\definecolor{stairs}{HTML}{F28A13}
\definecolor{embedding}{HTML}{A744D5}
\definecolor{guide}{HTML}{C84887}
\definecolor{purple}{HTML}{A744D5}
\definecolor{orange}{HTML}{E78A18}

\usepackage{graphics,graphicx}

\let\oldtocsection=\tocsection
\let\oldtocsubsection=\tocsubsection

\renewcommand{\tocsection}[2]{\hspace{0em}\oldtocsection{#1}{#2}}
\renewcommand{\tocsubsection}[2]{\hspace{1em}\oldtocsubsection{#1}{#2}}

\tikzset{node distance=3cm, auto}

\makeatletter
\def\@secnumfont{\bfseries}
\def\section{\@startsection{section}{1}%
  \z@{.7\linespacing\@plus\linespacing}{.5\linespacing}%
  {\normalfont\Large\bfseries}}
\def\subsection{\@startsection{subsection}{2}%
  \z@{.5\linespacing\@plus.7\linespacing}{-.5em}%
  {\normalfont\large\bfseries}}
  \def\subsubsection{\@startsection{subsubsection}{3}%
  \z@{.5\linespacing\@plus.7\linespacing}{-.5em}%
  {\normalfont\bfseries}}
\makeatother

\newtheorem{thm}{Theorem}[section]
\newtheorem{lemma}[thm]{Lemma}
\newtheorem{prop}[thm]{Proposition}
\newtheorem{cor}[thm]{Corollary}

\newtheorem{mainthm}{Theorem} 

\newtheorem{maincor}[mainthm]{Corollary}
\newtheorem{mainconj}[mainthm]{Conjecture}

\newtheorem{definition}[thm]{Definition}

\theoremstyle{remark}
\numberwithin{equation}{section} 
 \numberwithin{figure}{section}

\newtheoremstyle{customremark}
{3pt}
{3pt}
{}
{}
{\bfseries}
{.}
{.5em}
{}
\theoremstyle{customremark}
\newtheorem{rmk_no_diamond}[thm]{Remark}
\newenvironment{rmk}{\begin{rmk_no_diamond} } {\hfill$\er$ \end{rmk_no_diamond}}
\newtheorem{example_no_diamond}[thm]{Example}
\newenvironment{example}{\begin{example_no_diamond} } {\hfill$\er$ \end{example_no_diamond}}

\newcommand{\bF}{{\bf F}}

\newcommand{\Tt}{\mathcal{T}}

\newcommand{\bE}{{\bf{E}}}

\newcommand{\ov}{\overline}

\newcommand{\de}{{\delta}}

\newcommand{\la}{{\lambda}}

\newcommand{\p}{{\partial}}

\newcommand{\er}{{\Diamond}}
\newcommand{\Z}{\mathbb{Z}}

\newcommand{\R}{\mathbb{R}}

\newcommand{\C}{\mathbb{C}}
\newcommand{\CP}{\mathbb{CP}}

\newcommand{\eps}{\varepsilon}

\newcommand{\acc}{\mathrm{acc}}
\newcommand{\sembeds}{\stackrel{s}{\hookrightarrow}}

\newcommand{\vol}{\mathrm{Vol}}
\newcommand{\per}{\mathrm{Per}}
\makeatletter
\newcommand{\dashover}[2][\mathop]{#1{\mathpalette\df@over{{\dashfill}{#2}}}}
\newcommand{\fillover}[2][\mathop]{#1{\mathpalette\df@over{{\solidfill}{#2}}}}
\newcommand{\df@over}[2]{\df@@over#1#2}
\newcommand\df@@over[3]{%
  \vbox{
    \offinterlineskip
    \ialign{##\cr
      #2{#1}\cr
      \noalign{\kern1pt}
      $\m@th#1#3$\cr
    }
  }%
}
\newcommand{\dashfill}[1]{%
  \kern-.5pt
  \xleaders\hbox{\kern.5pt\vrule height.4pt width \dash@width{#1}\kern.5pt}\hfill
  \kern-.5pt
}
\newcommand{\dash@width}[1]{%
  \ifx#1\displaystyle
    2pt
  \else
    \ifx#1\textstyle
      1.5pt
    \else
      \ifx#1\scriptstyle
        1.25pt
      \else
        \ifx#1\scriptscriptstyle
          1pt
        \fi
      \fi
    \fi
  \fi
}
\newcommand{\solidfill}[1]{\leaders\hrule\hfill}
\makeatother

\date{\today}
\title[Full fillings and triangular limits of almost toric diagrams]{Full fillings at the accumulation point and triangular limits of almost toric diagrams}
\author{Nicki Magill}
\address{Mathematics Department, UC Berkeley}
\email{nmagill@berkeley.edu}
\thanks{NM thanks the NSF for their support under the agreement DMS-2402169}

\begin{document}
\begin{abstract}
Let $M$ be a closed symplectic four-manifold admitting an almost toric fibration over a disk. The ellipsoid embedding function of $M$ has an infinite staircase if it has infinitely many nonsmooth points. Cristofaro-Gardiner, Holm, Mandini, and Pires showed that these points must then accumulate at a solution of an explicit quadratic equation, and that at this accumulation point the ellipsoid embedding function agrees with the volume obstruction. We call this a full filling at the accumulation point, and study it in terms of almost toric fibrations. We show that if $M$ admits a sequence of almost toric base diagrams, each with a Delzant corner, converging nontrivially to a triangle, then $M$ has a full filling at its accumulation point. For the one-point blowups $H_b$ of $\mathbb{CP}^2$, we also prove the converse. As a consequence, $H_b$ has an infinite staircase if and only if it admits such a sequence converging to a triangle whose hypotenuse has irrational slope.
\end{abstract}

\maketitle

 \section{Introduction}

Let $(M,\omega)$ be a closed symplectic four-manifold admitting
an almost toric fibration over a disk. The classification of Leung and Symington implies that $M$ is diffeomorphic to
$\CP^2$, a blowup of $\CP^2$, or $S^2\times S^2$ \cite[Theorem 6.12]{LS}. Much work has gone into studying
the ellipsoid embedding function for such manifolds. This function
determines when the four-dimensional ellipsoid
\[
E(a,b)=
\left\{(z_1,z_2)\in\C^2 \ \middle|\
\pi\left(\frac{|z_1|^2}{a}+\frac{|z_2|^2}{b}\right)<1
\right\}
\]
can be symplectically embedded into $(M,\omega)$. The \textbf{ellipsoid embedding function} of $(M,\omega)$ is defined by
\[
c_M(z):=
\inf\{\lambda>0 \mid E(1,z)\sembeds \lambda M\},
\qquad z\geq1,
\]
where $\lambda M$ denotes $(M,\lambda\omega)$ and $\sembeds$ denotes a symplectic embedding. Throughout, we
suppress the symplectic form from the notation. We use the
normalization $\vol(M):=\int_M\omega^2.$

For a symplectic embedding to exist, the volume of the domain must be no greater than the volume of the target. This 
 gives the lower bound
\[
c_M(z)\geq V_M(z):=\sqrt{\frac{z}{\vol(M)}}.
\]

Work of Cristofaro-Gardiner, Holm, Mandini, and Pires shows the function $c_M(z)$ equals the volume obstruction for sufficiently
large $z$, and is locally piecewise linear wherever it is strictly
above the volume obstruction \cite[Proposition~2.1]{CGHMP}.
When its graph has infinitely many nonsmooth points accumulating
at a finite value, we say that $c_M$ has an \textbf{infinite staircase}.
The first example was found by McDuff and Schlenk in their study
of ellipsoid embeddings into a ball, equivalently into $\CP^2$
\cite{FibStairs}. Since then, infinite staircases have been studied
for a variety of symplectic targets \cite{REU,usher,CGMM,MMW,MM,CV,CGHMP,MPW,MarkovStairs}.

Here, we focus on the accumulation point theorem of
Cristofaro-Gardiner, Holm, Mandini, and Pires
\cite[Theorem 1.13 and Section 1.3]{CGHMP}. We define the perimeter of $M$ to be 
$c_1(M,\omega)\cdot[\omega].$
The theorem states that if $c_M$ has an infinite staircase,
then it accumulates at a real solution $z_\infty \geq 1$ of
\begin{equation}\label{eq:acc}
    z+\frac1z=\frac{\per(M)^2}{\vol(M)}-2.
\end{equation}

Furthermore,
\[
c_M(z_\infty)=V_M(z_\infty).
\]
Whenever the equation has a real solution greater than or equal to $1$,
we denote this solution by $\acc(M)$ and call it the
\textbf{accumulation point}, whether or not $c_M$ has an infinite
staircase. We refer to the equality
\[
c_M(\acc(M))=V_M(\acc(M))
\]
as \textbf{a full filling at the accumulation point}.

The accumulation point theorem has been very fruitful in proving that various $M$ do not have staircases, as ideally one only needs to compute the function at $\acc(M)$. For the one-point blowups
\[
H_b:=\CP^2_1\#\overline{\CP}{}^2_b,
\qquad b\in[0,1),
\]
with the convention $H_0=\CP^2_1$, a series of papers establishes a converse to the accumulation point
theorem \cite{MMW,MM,MPW,ICERM}.
Writing $c_b:=c_{H_b}$ and $\acc(b):=\acc(H_b)$, the function $c_b$
has an infinite staircase if and only if $\acc(b)$ is irrational
and
\[
c_b(\acc(b))=V_{H_b}(\acc(b)).
\]

In this paper, we study the geometric meaning of full filling
at the accumulation point. We show that if $M$ admits a nontrivial sequence
of almost toric base diagrams converging to a triangle, then $M$ has a full filling at the accumulation point; see Theorem~\ref{thm:general}. For $H_b$, we also prove the converse, giving
a characterization of a full filling at the accumulation point in terms of triangular limits of almost toric base diagrams; see Theorem~\ref{thm:Hb}.

An almost toric fibration in our setting is represented by an \textbf{almost toric base diagram} (ATBD), which is a convex polygon decorated
with nodes and cuts. Throughout, each cut is a straight segment along the eigenline of its node, joining the node to the boundary of the polygon. Interior points of the polygon other than the nodes are regular values of the
fibration, with Lagrangian torus fibers. Nodes correspond to focus-focus singularities, while boundary points correspond to
elliptic singularities. The cuts allow us to represent the integral
affine structure on the regular base in planar coordinates. Ignoring the decorations gives the \textbf{underlying polygon}
of an ATBD $Q$.  For an ATBD for $(M,\omega)$, the area of the polygon is $\tfrac12\vol(M)$, and the affine perimeter of the polygon is $\per(M).$

If $Q$ has a Delzant vertex, we say that $Q$ is in
\textbf{standard position} with respect to this vertex if the vertex is at the origin, no cut ends at it, and its adjacent edges lie along the positive coordinate axes.

Ellipsoid embeddings using these diagrams were first considered by Casals--Vianna \cite{CV} and Cristofaro-Gardiner--Holm--Mandini--Pires \cite{CGHMP}. The closure of the moment image of $E(a,b)$ is the triangle
\[
\Delta=\operatorname{conv}\{(0,0),(a,0),(0,b)\}.
\]
An ATBD in standard position whose axis edges end at $(a,0)$ and $(0,b)$ gives a symplectic embedding of $(1-\eta)E(a,b)$ into $M$ for every $0<\eta<1$ \cite[Proposition~2.35]{CGHMP}; we call these \emph{visible embeddings}. In particular, if a sequence of ATBDs in standard position
converges to $\Delta$, with $a\geq b>0$, then $M$ contains
ellipsoids of aspect ratio $a/b$ whose volumes approach that
of $M$. Our result identifies this aspect ratio with $\acc(M)$
when the convergence is nontrivial in the following sense.

\begin{definition}\label{def:converge}
Let $Q_k$ be a sequence of almost toric base diagrams in standard
position, and let $X_k,Y_k$ be the endpoints other than the origin
of their horizontal and vertical axis edges, respectively.
We say that $Q_k$ \emph{converges nontrivially} to
\[
\Delta=\operatorname{conv}\{(0,0),(A,0),(0,B)\},
\qquad A,B>0,
\]
if
\begin{enumerate}
    \item each underlying polygon occurs for only finitely many $k$;
    \item the underlying polygons converge to $\Delta$ in the
    Hausdorff topology;
    \item $X_k\to(A,0)$ and $Y_k\to(0,B)$.
\end{enumerate}
\end{definition}

Condition (1) rules out trivial examples. For instance, the standard triangle with vertices $(0,0)$, $(1,0)$, and $(0,1)$ is a toric diagram for $\CP^2$, and the constant sequence of this diagram satisfies (2) and (3), with $A/B=1$. When $A/B$ is irrational, condition (1) follows from (2). A polygon occurring infinitely often would have to equal the limiting triangle, whereas every edge of an almost toric base diagram has rational slope. In the rational case, condition (1) is essential.

\begin{mainthm}\label{thm:general}
Let $(M,\omega)$ be a closed symplectic four-manifold admitting
a sequence of almost toric base diagrams in standard position
that converges nontrivially to a triangle
\[
\Delta=\operatorname{conv}\{(0,0),(A,0),(0,B)\},
\qquad A,B>0.
\]
Then $\acc(M)$ is defined,
\[
\max\{A/B,B/A\}=\acc(M),
\]
and $M$ admits a full filling at its accumulation point:
\[
c_M(\acc(M))
=\sqrt{\frac{\acc(M)}{\vol(M)}}.
\]
\end{mainthm}

\begin{rmk}
The distinction between rational and irrational limiting slopes is important in the proof of Theorem~\ref{thm:general}. Visible embeddings give a full filling at the limiting aspect ratio in either case; identifying this ratio with the accumulation point requires showing that the total affine length of the non-axis edges tends to zero. When $A/B$ is irrational, this follows from an elementary argument using primitive integral edge directions. When $A/B$ is rational, an additional argument using visible symplectic surfaces is needed to rule out edges of affine length bounded away from zero. In the rational case, the affine perimeter is discontinuous at the limit. The diagrams all have affine perimeter \(\per(M)\), and the proof shows that \(\per(M)=A+B\). However, the limiting triangle has affine perimeter \(A+B+\ell_{\mathrm{aff}}(\text{hypotenuse})>A+B\), since its hypotenuse has rational slope.Example~\ref{ex:perimeterJump} illustrates this phenomenon for $H_{1/5}$.
\end{rmk}

For $H_b$, we also prove the converse.

\begin{mainthm}\label{thm:Hb}
Let $b\in[0,1)$. The following are equivalent:
\begin{enumerate}
    \item $H_b$ admits a full filling at its accumulation point:
    \[
    c_b(\acc(b))=\sqrt{\frac{\acc(b)}{1-b^2}}.
    \]
    \item There exists a sequence of almost toric base diagrams
    for $H_b$ in standard position that converges nontrivially
    to a triangle.
\end{enumerate}
\end{mainthm}

The implication (2)$\Rightarrow$(1) is a special case of Theorem~\ref{thm:general}. For (1)$\Rightarrow$(2), we use the classification of the parameters $b$ for which $H_b$ admits a full filling at its accumulation point: by \cite[Remark~1.2.2]{MPW}, building on \cite{MMW}, this happens exactly when $c_b$ has an infinite staircase or $b$ is one of the special rational values $b_i$ defined in \eqref{eq:special}. We then construct the required sequences case by case. For $b=0$ and $b=1/3$ they are given by Casals--Vianna \cite{CV} and Cristofaro-Gardiner--Holm--Mandini--Pires \cite{CGHMP}, and for ascending principal staircases by the author in \cite{M1,M2}. We treat the descending principal and nonprincipal staircases in Section~\ref{ss:descending}, and the special rational values in Section~\ref{ss:specialRational}.

\begin{maincor}\label{cor:Hb}
For $b\in[0,1)$, the function $c_b$ has an infinite staircase if and only if there exists a sequence of almost toric base diagrams for $H_b$ in standard position that converges nontrivially to a triangle whose hypotenuse has irrational slope.
\end{maincor}

Indeed, by Theorem~\ref{thm:general} the limiting triangle satisfies $\max\{A/B,B/A\}=\acc(b)$, so its hypotenuse has irrational slope exactly when $\acc(b)$ is irrational; Corollary~\ref{cor:Hb} then follows from Theorem~\ref{thm:Hb} and the characterization of infinite staircases recalled above.

The diagrams we construct for $H_b$ also record the staircase itself. For a staircase parameter $b$, the nodal ray at the vertex $Y$ on the $y$-axis of the $k$-th diagram in our sequence has slope $-p_k/q_k$, where the sequence $p_k/q_k$ contains subsequences giving the $z$-coordinates of staircase corners of $c_b$ \cite{MMW,MPW}. These converge to $\acc(b)$, the ratio of the side lengths of the limiting triangle. Thus the limiting triangle gives the full filling at $\acc(b)$, while the sequence of diagrams converging to it records the staircase. For the special rational values $b_i$, which have no staircase, the nodal rays instead record the steps in the special monotone staircase at $b=1/3$ and shifts of the basic blocking classes from \cite{MMW}. As we discuss in Remark~\ref{rmk:descendVisible}, unlike in the ascending case, we do not expect these mutation sequences to give visible obstructions in the sense of McDuff--Siegel \cite{McDuffSiegel2025}.

Motivated by Theorem~\ref{thm:Hb}, we conjecture that the converse
of Theorem~\ref{thm:general} holds more generally.

\begin{mainconj}\label{conj:general}
Let $(M,\omega)$ be a closed symplectic four-manifold admitting
an almost toric base diagram in standard position, and suppose
that $\acc(M)$ is defined. The following are equivalent:
\begin{enumerate}
    \item $M$ admits a full filling at its accumulation point:
    \[
    c_M(\acc(M))=V_M(\acc(M)).
    \]
    \item There exists a sequence of almost toric base diagrams
    for $(M,\omega)$ in standard position that converges
    nontrivially to a triangle.
\end{enumerate}
\end{mainconj}

Conjecture~\ref{conj:general} holds in all monotone cases with infinite staircases: the sequences of almost toric base diagrams constructed by Casals--Vianna \cite{CV} and Cristofaro-Gardiner--Holm--Mandini--Pires \cite{CGHMP} converge nontrivially to a triangle. It also holds for one example of $S^2_1 \times S^2_b$ for irrational $b$ studied in \cite{REU}. Moreover, the author has verified the implication (1)$\Rightarrow$(2) for a family of staircases in $\CP^2\#2\overline{\CP}^2$.
\par\medskip
\noindent\textbf{AI assistance.} I used ChatGPT and Claude to assist with computations leading to the sequences for the special rational $b$-values in Proposition~\ref{lem:rationalTriangle}, to identify and help resolve a gap in Lemma~\ref{lem:finite}, and to simplify computations and improve exposition. These tools also identified smaller errors, checked computations, and converted hand-drawn figures into TikZ code. I checked the resulting arguments and computations and take responsibility for the final manuscript.

\section{Triangle convergence}

In this section we prove Theorem~\ref{thm:general}. We first must define affine length: for an edge $e$ of rational slope with endpoints $R,S$, write
\[
S-R=\ell u,
\]
where $u\in\Z^2$ is primitive and $\ell>0$.
The \emph{affine length} of $e$ is
$\ell_{\mathrm{aff}}(e):=\ell$.
If an edge has irrational slope, we define its affine length to be zero.
We write $\per(Q)$ for the \emph{affine perimeter} of a polygon $Q$,
that is, the sum of the affine lengths of its edges.
If $Q$ is an almost toric base diagram for $(M,\omega)$, then
$\per(Q)=\per(M)$.

The idea of the theorem comes from the observation in \cite[Lemma 4.1.6]{M2}: if a sequence of diagrams converges to a triangle and the total affine length of the edges converging to its hypotenuse tends to zero, then the side lengths of the triangle are related by the accumulation point equation \eqref{eq:acc}. The first step of the proof below is the general form of this observation. We give the proof assuming Lemma~\ref{lem:tendtozero}, which shows that the total affine length of the non-axis edges tends to zero, and then prove Lemmas~\ref{lem:finite} and~\ref{lem:tendtozero}.

\begin{proof}[Proof of Theorem~\ref{thm:general}]
For each $k$, let $A_k:=|OX_k|$ and $B_k:=|OY_k|$. Twice the area of $Q_k$ equals $\vol(M)$, and its affine perimeter equals $\per(M)$. By Lemma~\ref{lem:tendtozero}, the total affine length of the non-axis edges of $Q_k$ tends to zero, so passing to the limit gives
\[
\vol(M)=AB,\qquad \per(M)=\lim_{k\to\infty}(A_k+B_k)=A+B.
\]
Hence
\[
\frac AB+\frac BA=\frac{(A+B)^2}{AB}-2=\frac{\per(M)^2}{\vol(M)}-2,
\]
so $A/B$ and $B/A$ solve \eqref{eq:acc}. In particular, $\acc(M)$ is defined and equals $\max\{A/B,B/A\}$. By symmetry we may assume $A\geq B$, so $\acc(M)=A/B$.
 By \cite[Proposition~2.35]{CGHMP}, for any $0<\eta<1,$ there exists a symplectic embedding of the ellipsoid $(1-\eta)E(A_k,B_k)$ into $M.$
 Set $z_k:=\max\{A_k/B_k,\,B_k/A_k\}\geq1$ and $m_k:=\min\{A_k,B_k\}$. Since $E(A_k,B_k)$ is symplectomorphic to $m_k\,E(1,z_k)$, we obtain $c_M(z_k)\leq \frac{1}{(1-\eta)m_k}$ for every $\eta$, and hence $c_M(z_k)\leq 1/m_k$.
 As $k\to\infty$, we have $z_k\to A/B=\acc(M)$ and $m_k\to B$. Since $c_M$ is continuous \cite[Proposition~2.1]{CGHMP}, it follows that
 \[
 c_M(\acc(M))\leq \frac1B=\sqrt{\frac{A/B}{AB}}=V_M(\acc(M)),
 \]
 as desired.
\end{proof}

It remains to prove that the total affine length of the non-axis edges
tends to zero. The idea is as follows. Suppose, for a contradiction,
that some non-axis edge has affine length bounded below along a
subsequence. Then we conclude that infinitely many diagrams have an
edge with a fixed normal vector.

For an edge $e$ of a convex polygon $Q$ with primitive outward
normal $(r,s)$, its \emph{supporting line} is the line
\[
rx+sy=h
\]
containing $e$, with $Q$ contained in the half-plane
$rx+sy\leq h$.
We will show that, for a fixed normal $(r,s)$ with $r,s>0$, only finitely
many values of $h$ in a bounded interval can occur among the diagrams
under consideration.
If the limiting hypotenuse has irrational slope, boundedness of the primitive edge directions already gives a contradiction. In the rational case, we additionally use the finiteness of the possible supporting heights, proved below using visible symplectic surfaces. See Example~\ref{ex:perimeterJump} for a rational-slope limit.

We recall the construction of visible symplectic surfaces, following
\cite[Section 7]{S} and \cite[Section 5.1a]{McDuffSiegel2025}. A \textbf{covector-decorated path}
$(\gamma,\alpha)$ consists of an embedded path $\gamma\colon[0,1]\to Q$
whose interior lies in the regular part of the base, together with a
primitive integral covector $\alpha$ along $\gamma$ that is constant in
integral affine coordinates. In action-angle coordinates away from the
nodal cuts, the covector $\alpha = a_1\,dx + a_2\,dy$ corresponds to
the vector field
\[
  X_\alpha = a_1\frac{\partial}{\partial\theta_1}
           + a_2\frac{\partial}{\partial\theta_2},
\]
where $\theta_1,\theta_2$ are period-one angular coordinates on the torus
fibers. The orbits of $X_\alpha$ are embedded circles in the fibers, and
choosing one such circle in each fiber over the interior of $\gamma$,
smoothly in $t$, gives a \textbf{visible cylinder} $C$.
If $\alpha(\gamma'(t)) > 0$ for all $t\in(0,1)$, then $C$ is symplectic,
with area
\begin{equation}\label{eq:visibleArea}
  \omega(C) = \int_0^1 \alpha(\gamma'(t))\,dt,
\end{equation}
which equals $\alpha(\gamma(1)-\gamma(0))$ when $\gamma$ crosses no cuts.

Here, we only consider endpoints at a Delzant vertex or at an interior
point of an edge. At an interior point of an edge $e$ with primitive
outward normal $(r,s)$, the circle that collapses is the orbit of
$X_\alpha$ for $\alpha = r\,dx + s\,dy$. Hence if $\alpha$ is this
covector and $\gamma$ is a straight segment in the direction $(r,s)$
near its endpoint, the cylinder caps off smoothly to a symplectic disk.
At a Delzant vertex, choose integral affine coordinates in which the
vertex is the origin and the adjacent edges lie on the positive
coordinate axes. If $\alpha = r\,dx + s\,dy$ with $r,s>0$ and
$\gcd(r,s)=1$, the cylinder is completed by a symplectic disk with an
$(r,s)$-cusp at the vertex
\cite[Section 5.1a, items (ii) and (iv)]{McDuffSiegel2025}. We call the resulting
surface the \textbf{visible surface} associated to $(\gamma,\alpha)$.

\begin{lemma}\label{lem:finite}
Let $Q$ be an almost toric base diagram for $(M,\omega)$ in standard position, and let $e$ be an edge of $Q$ with primitive outward normal $(r,s)$, where $r,s>0$, lying on the line $rx+sy=h$. Then
\[
h\in\mathcal{A}_{rs}:=\{\omega(A)\mid A\in H_2(M;\Z),\ A^2=rs,\ \omega(A)>0\}.
\]
Moreover, for every $h_0>0$, the set $\mathcal{A}_{rs}\cap(0,h_0]$ is finite.
\end{lemma}
\begin{proof}
Let $\alpha=r\,dx+s\,dy$. Choose a point $P$ in the interior
of $e$ that is not an endpoint of a cut. For sufficiently
small $\de>0$, the point
\[
P_\de:=P-\de(r,s)
\]
lies in $\operatorname{int}(Q)$ and satisfies
$h-\de(r^2+s^2)>0$. Consider the broken line from $O$
to $P_\de$ and then to $P$. By convexity, its interior
lies in $\operatorname{int}(Q)$, and $\alpha$ evaluates
positively on the tangent vector of each segment.
Smoothing the corner at $P_\de$ gives an embedded path
$\gamma\colon[0,1]\to Q$ such that
$\alpha(\gamma'(t))>0$ for $t\in(0,1)$. We leave the path
unchanged near its endpoints, so it meets $e$ orthogonally
in the direction $(r,s)$.

The path meets $\partial Q$ only at $O$ and $P$, neither
of which is an endpoint of a cut. We may therefore perform
nodal slides \cite[Definition~6.1 and Proposition~6.2]{S} to shorten the cuts sufficiently that
$\gamma$ avoids all cuts and marked points. Equip
$\gamma$ with the constant covector $\alpha$.

Let $C$ be the visible surface corresponding to $(\gamma,\alpha)$.
By \cite[Section 5.1a, items (ii) and (iv)]{McDuffSiegel2025}, $C$ is a rational
symplectic curve with an $(r,s)$-cusp at the origin, and by
\eqref{eq:visibleArea}, $\omega(C)=h>0$.
By \cite[Proposition~8.2]{S}, the preimage of the boundary represents the Poincar\'e dual of $c_1(M)$. The intersections at the origin contribute $r+s$, and the cap at $e$ contributes one. Hence $c_1(C)=r+s+1$.
Since $C$ is rational and embedded away from its $(r,s)$-cusp,
the adjunction formula gives
\[
[C]^2=c_1(C)-2+(r-1)(s-1)=rs.
\]
Hence $h=\omega([C])\in\mathcal{A}_{rs}$.

For the second claim, fix $h_0>0$. It suffices to show that
\[
  \{A \in H_2(M;\Z) \mid 0 < \omega(A) \leq h_0,\ A^2 = rs\}
\]
is finite.

Let $\Omega = \operatorname{PD}([\omega])$, so that
$\Omega^2 = \int_M \omega^2 > 0$. Since $M$ admits an almost toric fibration over a disk, $b_2^+(M) = 1$ by the classification in \cite[Theorem~6.12]{LS}, so the intersection
form is negative definite on $\Omega^\perp \subset H_2(M;\R)$.
Writing
\[
  A = \frac{\omega(A)}{\Omega^2}\,\Omega + A',
  \qquad A' \in \Omega^\perp,
\]
we obtain
\[
  -A'^2 = \frac{\omega(A)^2}{\Omega^2} - rs
  \leq \frac{h_0^2}{\Omega^2} - rs.
\]
Hence the coefficient of $\Omega$ is bounded by $h_0/\Omega^2$, and since
$-A'^2$ is a positive definite quadratic form on $\Omega^\perp$, the
component $A'$ is also bounded. Thus all such classes lie in a compact
subset of $H_2(M;\R)$.
Since $H_2(M;\Z)$ is a lattice in $H_2(M;\R)$, this set
is finite, and hence so is its image under $\omega$.
\end{proof}

The following lemma shows that the total affine length of the
non-axis edges tends to zero.

\begin{lemma}\label{lem:tendtozero}
Let $Q_k\subset\R_{\geq0}^2$ be a sequence of ATBDs for a fixed
symplectic manifold $(M,\omega)$, standardized with a Delzant corner
at the origin and its adjacent edges along the coordinate axes.
Assume that $Q_k$ converges nontrivially to
\[
T=\operatorname{conv}\{(0,0),(A,0),(0,B)\},
\qquad A,B>0,
\]
in the sense of Definition~\ref{def:converge}.
Then the total affine length of the non-axis edges of $Q_k$
tends to zero.
\end{lemma}
\begin{proof}
The number of edges of each $Q_k$ is bounded independently of $k$ by $\chi(M)$ \cite[Section 2]{LS}. Thus it suffices to show that the maximum affine length of a non-axis edge tends to zero.
Suppose, for a contradiction, that along a subsequence there is a
non-axis edge with affine length $\ell_k\geq\de>0$.
Let its endpoints be $R_k,S_k$, and write
\[
R_k-S_k=\ell_k u_k,
\]
where $u_k$ is the primitive integral direction.
Since the diagrams converge to a triangle, their Euclidean edge
lengths are uniformly bounded. Thus
\[
\|u_k\|=\frac{\|R_k-S_k\|}{\ell_k}
\leq \frac{\|R_k-S_k\|}{\de}
\]
is bounded. Since $u_k$ is integral, there are only finitely many possibilities for it, so after passing to a subsequence we may assume it is constant and the corresponding primitive outward
normal $(r,s)$ is fixed.

The chosen edges therefore lie on supporting lines
\[
rx+sy=h_k
\]
for some $h_k$, and convexity implies that $Q_k$ lies in the
half-plane $rx+sy\leq h_k$.

Passing to a further subsequence, assume that $R_k\to R$ and
$S_k\to S$. The lower bound on $\ell_k$ implies that $R\ne S$.
Moreover, $h_k\to h_\infty$, and the line $rx+sy=h_\infty$
contains the segment $\overline{RS}$, which lies on the boundary of the limiting triangle.

This segment cannot lie on a coordinate axis. Indeed, if it lay
on the horizontal axis, the fixed outward normal would be $(0,-1)$.
The corresponding edges would then lie on lines $y=c_k$ with
$Q_k\subseteq\{y\geq c_k\}$. Since
$O\in Q_k\subseteq\mathbb R_{\geq0}^2$, we would have $c_k=0$,
contradicting the choice of a non-axis edge.
The vertical-axis case is identical.

Thus $\overline{RS}$ lies on the hypotenuse, and $r,s>0$.
Evaluating $rx+sy=h_\infty$ at the endpoints $(A,0)$ and $(0,B)$ of the hypotenuse gives $rA=h_\infty=sB$, so $A/B=s/r$ is rational. Hence if $A/B$ is irrational, we already have a contradiction. It remains to treat the case where $A/B$ is rational.
By Lemma~\ref{lem:finite}, each $h_k$ lies in the set $\mathcal{A}_{rs}$, which has only finitely many elements in a bounded interval, so we may assume $h_k=h$ is constant.
Hence the limiting triangle $T$ is the intersection of the
positive quadrant and $rx+sy\leq h$.

By convexity $Q_k\subseteq T$. Since all $Q_k$ represent the same
symplectic manifold, they have the same area, which also equals
the area of $T$. Thus $Q_k=T$, contradicting the assumption that
no underlying polygon occurs infinitely often.

This contradiction shows that the maximum affine length among the non-axis edges of $Q_k$ tends to zero. By the uniform bound on the number of edges, their total affine length also tends to zero.
\end{proof}

\begin{example}\label{ex:perimeterJump}
The rational case in the proof of Lemma~\ref{lem:tendtozero} does occur, for the special rational values $b_i$ of \eqref{eq:special}. In these cases, the affine perimeter of the limiting triangle is strictly larger than $\per(M)$. By the proof of Theorem~\ref{thm:general}, $\per(M)=A+B$, whereas the triangle $T$ has affine perimeter $A+B+\ell_{\mathrm{aff}}(\text{hypotenuse})$. If $A/B$ is irrational, the hypotenuse has irrational slope and affine length zero, so these agree; if $A/B$ is rational, the affine perimeter of $T$ is strictly larger than $\per(M)$.

Consider the special rational value $b=1/5$. In Proposition~\ref{lem:rationalTriangle}, we construct a sequence of diagrams $Q_k=OX_kV_kY_k$ where we let $A_k=|OX_k|$ and $B_k=|OY_k|.$ The limiting triangle has $A=12/5$ and $B=2/5$, so \[A+B=14/5=3-1/5=\per(H_{1/5}),\] while its hypotenuse $(-12/5,2/5)=\tfrac25(-6,1)$ has affine length $2/5$, giving the triangle affine perimeter $16/5$. The formulas in the proposition give
\[
A_k=\frac{12}{5}-\frac{1}{5(2k+2)},\qquad B_k=\frac25,\qquad \ell_{\mathrm{aff}}(V_kY_k)=\frac{1}{5(2k+4)},
\]
so $\ell_{\mathrm{aff}}(X_kV_k)=\per(H_{1/5})-A_k-B_k-\ell_{\mathrm{aff}}(V_kY_k)=\frac{1}{10(k+1)(k+2)}\to0$, even though the edges $X_kV_k$ converge to the entire hypotenuse. This is because affine length is not continuous under Hausdorff convergence: the directions of $X_kV_k$ approach $(-6,1)$ without ever equaling it, so their primitive direction vectors become arbitrarily long.
\end{example}

 \section{The sequences for $H_b$}

 In this section, we prove that (1) implies (2) in
Theorem~\ref{thm:Hb}. For each $b$ satisfying (1), we
construct an explicit sequence of almost toric base
diagrams that converges nontrivially to a triangle.

As explained after Theorem~\ref{thm:Hb}, by \cite[Remark~1.2.2]{MPW} it suffices to treat the values of $b$ for which $c_b$ has an infinite staircase and the special rational values, which we now define.

Let
\[
(\mathfrak{p}_0,\mathfrak{q}_0)=(6,1), \qquad
(\mathfrak{p}_{i+1},\mathfrak{q}_{i+1})=(6\mathfrak{p}_i-\mathfrak{q}_i,\mathfrak{p}_i),
\]
and set $\eps_i=(-1)^{i+1}$. The
\emph{special rational values} are
\begin{equation} \label{eq:special}
b_i=\frac{\mathfrak{p}_i+\mathfrak{q}_i+3\eps_i}{3(\mathfrak{p}_i+\mathfrak{q}_i)+\eps_i},
\qquad i\geq0.
\end{equation}
The first values are
\[
b_0=\frac15,\qquad
b_1=\frac{11}{31},\qquad
b_2=\frac{59}{179},\qquad
b_3=\frac{349}{1045}.
\]
They satisfy $\acc(b_i)=\mathfrak{p}_i/\mathfrak{q}_i$. Although $H_{b_i}$
admits a full filling at its accumulation point,
$c_{b_i}$ has no infinite staircase
\cite[Theorem~1.2.1]{MPW}.

We call $b$ a \emph{staircase parameter} if $c_b$ has an
infinite staircase. A staircase is called \emph{ascending}
if its infinitely many nonsmooth points approach the
accumulation point from the left, and \emph{descending}
if they approach it from the right.

Apart from $b=0$ and $b=1/3$, the staircase parameters
fall into two families: those with \emph{principal}
staircases, for which $c_b$ has an ascending staircase
or a descending staircase, but not both, and those with
\emph{nonprincipal} staircases, for which $c_b$ has both.

Since the cases $b=0$ and $b=1/3$ are treated by Casals--Vianna \cite{CV} and Cristofaro-Gardiner--Holm--Mandini--Pires \cite{CGHMP}, and ascending principal staircases by the author \cite{M1,M2}, it remains to consider descending principal staircases, nonprincipal staircases, and the special rational $b$-values.

We follow the methods of \cite{M2,M1}, constructing
sequences of diagrams by almost toric mutations.
The mutation sequences are guided by the recursive
construction of the staircases. We first recall
the necessary background on almost toric mutations
in Section~\ref{ss:mutationBackground} and on
staircases in Section~\ref{ss:stair}.
In Section~\ref{ss:descending}, we construct the
sequences for descending principal and nonprincipal
staircases. Finally, in
Section~\ref{ss:specialRational}, we construct the
sequences for the special rational $b$-values.

\begin{proof}[Proof of Theorem~\ref{thm:Hb}]
The implication $(2)\Rightarrow(1)$ follows from
Theorem~\ref{thm:general}.

For $(1)\Rightarrow(2)$, suppose that $H_b$ admits a
full filling at $\acc(b)$. By
\cite[Remark~1.2.2]{MPW}, either $c_b$ has an infinite
staircase or $b$ is a special rational value.

For $b=0$ and $b=1/3$, the required sequences are
given in \cite[Theorem~1.4]{CV} and \cite[Proposition~5.3]{CGHMP}. For the remaining staircase
parameters, the staircase is either ascending
principal, descending principal, or nonprincipal.
The ascending principal case follows from
\cite[Theorem~1.1]{M1} and \cite[Theorem~D]{M2}. In the descending principal and
nonprincipal cases, this is shown in Proposition~\ref{prop:descending}. 

Finally, if $b$ is a special rational value,
Proposition~\ref{lem:rationalTriangle} gives a sequence
of diagrams converging nontrivially to a triangle.
This proves $(2)$ in every case.
\end{proof}

 \subsection{Background on mutations for $H_b$} \label{ss:mutationBackground}
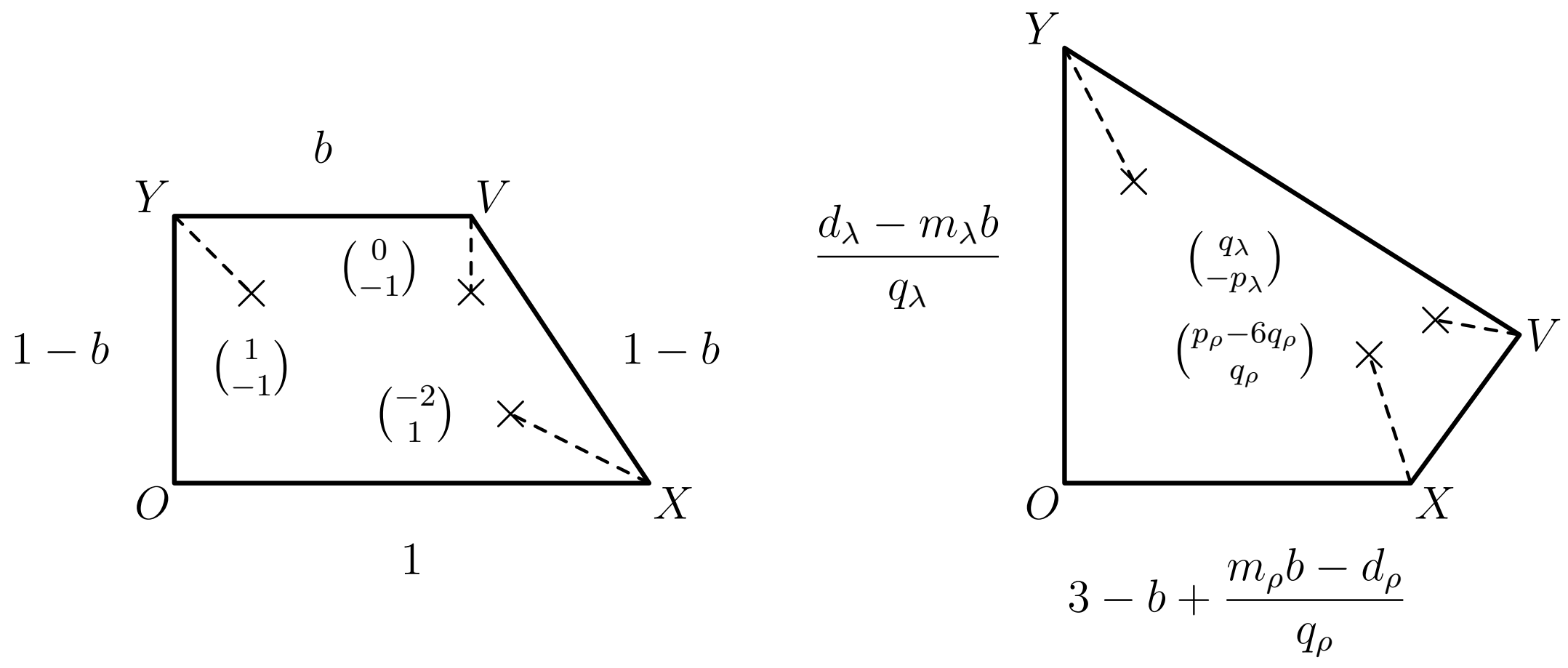
\begin{figure}
    \centering
\begin{tikzpicture}[scale=0.9,x=1cm,y=1cm,line cap=round,line join=round,
 every node/.style={font=\large},
 boundary/.style={line width=1.2pt},
 ray/.style={dashed,line width=.85pt},
 cross/.style={font=\Large,inner sep=0pt},
 vertex/.style={inner sep=0pt,outer sep=0pt}]

\begin{scope}
\coordinate (O) at (0,0);
\coordinate (X) at (4.8,0);
\coordinate (V) at (3,2.7);
\coordinate (Y) at (0,2.7);

\draw[boundary] (O)--(X)--(V)--(Y)--cycle;

\node[vertex,below left=1pt] at (O) {$O$};
\node[vertex,below right=1pt] at (X) {$X$};
\node[vertex,above right=1pt] at (V) {$V$};
\node[vertex,above left=1pt] at (Y) {$Y$};

\node[left=13pt] at ($(O)!.5!(Y)$) {$1-b$};
\node[above=10pt] at ($(Y)!.5!(V)$) {$b$};
\node[right=12pt] at ($(V)!.5!(X)$) {$1-b$};
\node[below=12pt] at ($(O)!.5!(X)$) {$1$};

\draw[ray] (Y)--(.78,1.92);
\node[cross] at (.78,1.92) {$\times$};
\node[below=8pt] at (.78,1.92) {$\binom{1}{-1}$};

\draw[ray] (V)--(3,1.93);
\node[cross] at (3,1.93) {$\times$};
\node[left=9pt] at (3,2.16) {$\binom{0}{-1}$};

\draw[ray] (X)--(3.4,.7);
\node[cross] at (3.4,.7) {$\times$};
\node[left=10pt] at (3.4,.7) {$\binom{-2}{1}$};
\end{scope}

\begin{scope}[xshift=9cm]
\coordinate (O) at (0,0);
\coordinate (X) at (3.5,0);
\coordinate (V) at (4.6,1.5);
\coordinate (Y) at (0,4.4);

\draw[boundary] (O)--(X)--(V)--(Y)--cycle;

\node[vertex,below left=1pt] at (O) {$O$};
\node[vertex,below right=1pt] at (X) {$X$};
\node[vertex,right=1.5pt] at (V) {$V$};
\node[vertex,above left=1pt] at (Y) {$Y$};

\node[left=12pt] at ($(O)!.52!(Y)$)
 {$\dfrac{d_\lambda-m_\lambda b}{q_\lambda}$};
\node[below=13pt] at ($(O)!.5!(X)$)
 {$3-b+\dfrac{m_\rho b-d_\rho}{q_\rho}$};

\draw[ray] (Y)--(.7,3.05);
\node[cross] at (.7,3.05) {$\times$};
\node[below right=9pt] at (.7,3.05)
 {$\binom{q_\lambda}{-p_\lambda}$};

\draw[ray] (X)--(3.08,1.3);
\node[cross] at (3.08,1.3) {$\times$};
\node[left=9pt] at (3.08,1.3)
 {$\binom{p_\rho-6q_\rho}{q_\rho}$};

\draw[ray] (V)--(3.75,1.65);
\node[cross] at (3.75,1.65) {$\times$};
\end{scope}

\end{tikzpicture}
    \caption{Left: the starting diagram $Q_b^0$, obtained by
performing three nodal trades on the moment polygon of $H_b$,
with its nodal rays and affine edge lengths labeled.
Right: a quadrilateral $Q_b(\Tt)$ associated with a triple
$\Tt=((p_\la,q_\la,t_\la),(p_\mu,q_\mu,t_\mu),
(p_\rho,q_\rho,t_\rho))$, showing how the triple determines
the affine lengths $|OX|$ and $|OY|$ and the nodal rays
at $X$ and $Y$. The quantities $(d,m)$ are defined from
$(p,q,t)$ in \eqref{eq:degree}.
The starting diagram corresponds to
$\Tt_0=((1,1,2),(2,1,-1),(4,1,1))$, for which
$(d_\la,m_\la)=(1,1)$ and $(d_\rho,m_\rho)=(2,1)$.
The triple also determines the nodal ray at $V$ and
the position of $V$; see
\cite[Definition~2.3.1 and Lemma~2.3.3]{M2}.}

    \label{fig:starting}
\end{figure}

 For background on almost toric fibrations, we refer the reader to Symington~\cite{S}, Leung--Symington~\cite{LS}, and the book by Evans~\cite{Evans}. We recall the almost toric base diagrams for $H_b$
and their mutations from \cite{M2}. The initial diagram $Q_b^0$ is shown in
Figure~\ref{fig:starting}; it is in standard position. 
Write
\[
Q=OXVY,
\qquad
X=(A,0),
\qquad
Y=(0,B).
\]
Thus $A=|OX|$ and $B=|OY|$. The area of the
diagram is $(1-b^2)/2$, and its affine perimeter is $\per(H_b)=3-b$.

Each diagram we consider has three nodal rays,
based at the vertices $X,Y,V$. We can mutate at each of these nodes, in the following sense.

\begin{definition}\label{def:mutation}
Let $Q$ be an almost toric base diagram in standard position, and let
$\mathfrak{b}$ be a node of $Q$ whose nodal ray has direction $(q,-p)$,
with $p,q$ coprime. The associated \emph{mutation matrix} is
\[
  M_{p,q} =
  \begin{pmatrix}
    1-pq & -q^2\\
    p^2 & 1+pq
  \end{pmatrix}
  \in \mathrm{SL}(2,\Z),
  \qquad\text{i.e.}\qquad
  M_{p,q}\,u = u - \det\bigl((q,-p),u\bigr)\,(q,-p),
\]
so that $M_{p,q}$ fixes $(q,-p)$.

Let $L$ be the eigenline through the node $\mathfrak b$,
oriented by the primitive vector $v=(q,-p)$ pointing
from the boundary vertex toward the node. Assume that
$O\notin L$, and let $\mathcal{H}$ be the closed half-plane
bounded by $L$ that does not contain $O$.

Set
\[
M_v=
\begin{cases}
M_{p,q}, & \text{if $\mathcal{H}$ lies to the left of the oriented line $L$},\\
M_{p,q}^{-1}, & \text{if $\mathcal{H}$ lies to the right of the oriented line $L$}.
\end{cases}
\]
The mutation at $\mathfrak b$ is obtained by reversing
the cut at $\mathfrak b$ and applying the affine shear
\[
P\longmapsto \mathfrak b+M_v(P-\mathfrak b)
\]
to $Q\cap \mathcal{H}$, leaving the other half unchanged.
Nodes and cuts in $\mathcal{H}$ are carried along by this shear.
Since the shear fixes $L$ pointwise, the two pieces
agree along $L$.
\end{definition}

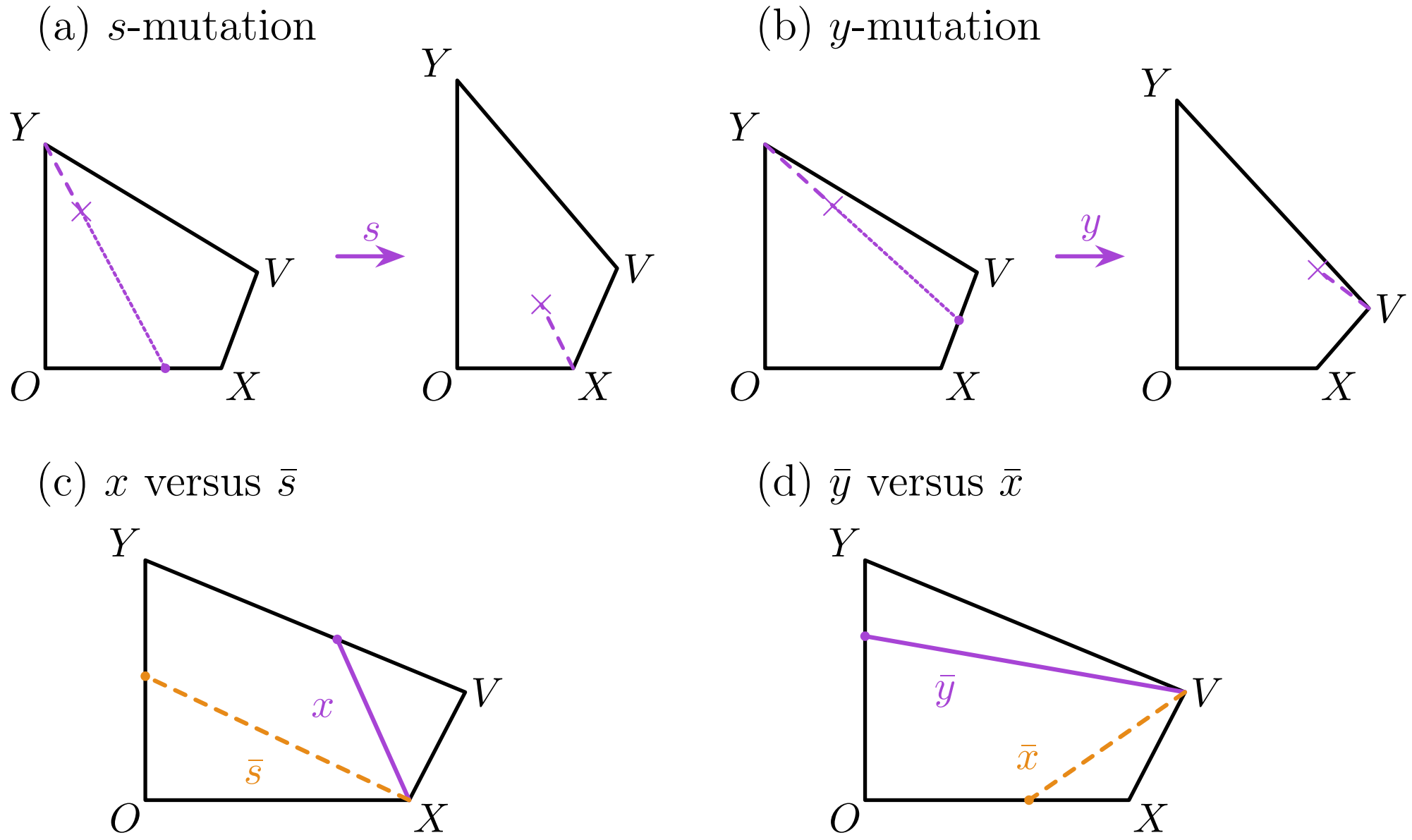
\begin{figure}
    \centering
\begin{tikzpicture}[scale=0.75,x=1cm,y=1cm,line cap=round,line join=round,>=Stealth,
 boundary/.style={line width=1pt},cut/.style={purple,line width=1pt,dashed},
 extension/.style={purple,line width=.85pt,densely dotted},
 every node/.style={font=\normalsize}]

\newcommand{\vertices}[4]{
 \node[below left=1pt,inner sep=0pt,outer sep=0pt] at #1 {$O$};
 \node[below right=1pt,inner sep=0pt,outer sep=0pt] at #2 {$X$};
 \node[right=1.5pt,inner sep=0pt,outer sep=0pt] at #3 {$V$};
 \node[above left=1pt,inner sep=0pt,outer sep=0pt] at #4 {$Y$};}
\newcommand{\crossat}[1]{\node[purple,inner sep=0pt] at #1 {$\times$};}

\begin{scope}
\node[anchor=west,font=\large] at (-.3,4.25) {(a) $s$-mutation};
\draw[boundary] (0,0)--(2.2,0)--(2.65,1.2)--(0,2.8)--cycle;
\vertices{(0,0)}{(2.2,0)}{(2.65,1.2)}{(0,2.8)}
\draw[cut] (0,2.8)--(.45,1.96); \crossat{(.45,1.96)}
\draw[extension] (.45,1.96)--(1.5,0);
\fill[purple] (1.5,0) circle (1.7pt);
\draw[purple,->,line width=1.1pt] (3.65,1.4)--(4.5,1.4) node[midway,above] {$s$};
\begin{scope}[xshift=.85cm]
\draw[boundary] (4.3,0)--(5.75,0)--(6.3,1.25)--(4.3,3.6)--cycle;
\vertices{(4.3,0)}{(5.75,0)}{(6.3,1.25)}{(4.3,3.6)}
\draw[cut] (5.75,0)--(5.35,.8); \crossat{(5.35,.8)}
\end{scope}
\end{scope}

\begin{scope}[xshift=9cm]
\node[anchor=west,font=\large] at (-.3,4.25) {(b) $y$-mutation};
\draw[boundary] (0,0)--(2.2,0)--(2.65,1.2)--(0,2.8)--cycle;
\vertices{(0,0)}{(2.2,0)}{(2.65,1.2)}{(0,2.8)}
\draw[cut] (0,2.8)--(.85,2.03); \crossat{(.85,2.03)}
\draw[extension] (.85,2.03)--(2.425,.6);
\fill[purple] (2.425,.6) circle (1.7pt);
\draw[purple,->,line width=1.1pt] (3.65,1.4)--(4.5,1.4) node[midway,above] {$y$};
\begin{scope}[xshift=.85cm]
\draw[boundary] (4.3,0)--(6.05,0)--(6.7,.75)--(4.3,3.35)--cycle;
\vertices{(4.3,0)}{(6.05,0)}{(6.7,.75)}{(4.3,3.35)}
\draw[cut] (6.7,.75)--(6.06,1.23); \crossat{(6.06,1.23)}
\end{scope}
\end{scope}

\begin{scope}[xshift=1.25cm,yshift=-5.4cm]
\node[anchor=west,font=\large] at (-1.55,3.95) {(c) $x$ versus $\bar s$};
\draw[boundary] (0,0)--(3.3,0)--(4,1.35)--(0,3)--cycle;
\vertices{(0,0)}{(3.3,0)}{(4,1.35)}{(0,3)}
\draw[purple,line width=1.15pt] (3.3,0)--(2.4,2.01);
\fill[purple] (2.4,2.01) circle (1.7pt);
\node[purple,left=5pt] at (2.78,1.15) {$x$};
\draw[orange,dashed,line width=1.15pt] (3.3,0)--(0,1.55);
\fill[orange] (0,1.55) circle (1.7pt);
\node[orange,below=5pt] at (1.35,.91) {$\bar s$};
\end{scope}

\begin{scope}[xshift=10.25cm,yshift=-5.4cm]
\node[anchor=west,font=\large] at (-1.55,3.95) {(d) $\bar y$ versus $\bar x$};
\draw[boundary] (0,0)--(3.3,0)--(4,1.35)--(0,3)--cycle;
\vertices{(0,0)}{(3.3,0)}{(4,1.35)}{(0,3)}
\draw[purple,line width=1.15pt] (4,1.35)--(0,2.05);
\fill[purple] (0,2.05) circle (1.7pt);
\node[purple,below=3pt] at (1.0,1.875) {$\bar y$};
\draw[orange,dashed,line width=1.15pt] (4,1.35)--(2.05,0);
\fill[orange] (2.05,0) circle (1.7pt);
\node[orange,below left=3pt] at (2.5,1) {$\bar x$};
\end{scope}
\end{tikzpicture}
    \caption{Panels (a) and (b) show the two possible mutations at $Y$ with their labeling of $s$ or $y$. Panels (c) and (d) label the possible nodal rays from $X$
and $V$, respectively, by the corresponding mutation,
according to which edge each ray intersects. For any given diagram, at most one of the two mutations
at each vertex is defined.}
    \label{fig:HbMutations}
\end{figure}

We call the mutations of Definition~\ref{def:mutation} \emph{geometric mutations}, to distinguish them from the \emph{algebraic mutations} of triples defined below; when the meaning is clear from context, we drop the adjective.
The directions of
the nodal rays and the side lengths of the quadrilateral
determine which edges the rays meet. We distinguish
the geometric mutations according to both the
vertex where the mutation is performed and the
edge met by its nodal ray.
We denote the mutation at $X$ by $x$ when its
nodal ray meets the interior of $\overline{VY}$,
and by $\ov{s}$ when it meets the interior of
$\overline{OY}$. Similarly, we denote the mutation
at $Y$ by $y$ when its nodal ray meets the interior
of $\overline{XV}$, and by $s$ when it meets the
interior of $\overline{OX}$. Finally, we denote
the mutation at $V$ by $\ov{x}$ when its nodal
ray meets the interior of $\overline{OX}$, and
by $\ov{y}$ when it meets the interior of
$\overline{OY}$. See Figure~\ref{fig:HbMutations}.
Thus, whether a particular geometric mutation is
well-defined depends
on the parameter $b$.

After each mutation, we standardize and relabel
the resulting diagram as $OXVY$. With this
convention, $\ov{x}$, $\ov{y}$, and $\ov{s}$
are the inverses of $x$, $y$, and $s$, respectively.

We compose mutation words from right to left.
For a finite word $w$, the notation $wQ_b^0$
means that we apply its geometric mutations in
this order. We say that $wQ_b^0$ is well-defined
if every mutation in the word is well-defined.
For example, $\ov{y}Q_b^0$ is not well-defined
for any $b\in(0,1)$, since the nodal ray at $V$
always meets the interior of $\overline{OX}$.

\medskip

Following \cite{M2}, we associate a triple
$\Tt_0\in(\Z^3)^3$ to $Q_b^0$. Together with $b$, this triple determines the diagram exactly: the side lengths and edge directions of the quadrilateral, and the directions of its nodal rays; see Figure~\ref{fig:starting}. After performing
a sequence of well-defined geometric mutations $w$
on $Q_b^0$, the affine lengths and directions of
the edges, as well as the directions of the nodal
rays, can be recovered from the triple $w\Tt_0$,
obtained by applying the corresponding algebraic
mutations to $\Tt_0$.

We call the triples obtained from $\Tt_0$ by the
algebraic mutations defined below \textbf{recursive
triples}. We now describe these mutations and
recall some properties of recursive triples
from \cite{M2}. Write
\[
\Tt=(\bE_\la,\bE_\mu,\bE_\rho),
\qquad
\bE_\alpha=(p_\alpha,q_\alpha,t_\alpha),
\quad \alpha\in\{\la,\mu,\rho\},
\]
and let $\eps\in\{1,-1\}$ be the sign of the
triple. We call $p_\alpha/q_\alpha$ the \emph{center} of $\bE_\alpha$. The initial triple is
\[
\Tt_0=((1,1,2),(2,1,-1),(4,1,1)),
\]
with sign $\eps=1$. For each entry $(p,q,t)$, define its \emph{degree coordinates}
\begin{equation} \label{eq:degree}
d=\frac{3(p+q)+t}{8},
\qquad
m=\frac{p+q+3t}{8}.
\end{equation}
By \cite[Definition~2.2.5 and Proposition~2.2.8]{M2}, for every entry of a recursive triple, $d$ and $m$ are integers, and
\[
3d-m=p+q,
\qquad
t=3m-d,
\qquad
d^2-m^2=pq-1.
\]

The corresponding algebraic mutations are
\[
x\Tt
=
(\bE_\la,\eps t_\la\bE_\mu-\bE_\rho,\bE_\mu),
\qquad
y\Tt
=
(\bE_\mu,\eps t_\rho\bE_\mu-\bE_\la,\bE_\rho),
\]
and we write $\ov{x}$ and $\ov{y}$
for the inverses of $x$ and $y$, e.g.,
$\ov{x}x\Tt=\Tt$.
All four mutations $x,y,\ov{x},\ov{y}$
preserve the sign $\eps$.
To describe $s$ and its inverse, set
\begin{equation}\label{eq:S}
S(p,q,t)=(6p-q,p,-t),
\qquad
S^{-1}(p,q,t)=(q,6q-p,-t).
\end{equation}
Then
\[
s\Tt=(\bE_\mu,\bE_\rho,S\bE_\la),
\qquad
\ov{s}\Tt=(S^{-1}\bE_\rho,\bE_\la,\bE_\mu).
\]
These mutations change the sign from $\eps$
to $-\eps$.

\begin{definition}[{\cite[Lemma~2.1.4 and Definition~2.3.1]{M2}}]\label{def:Qb}
For a recursive triple $\Tt=(\bE_\la,\bE_\mu,\bE_\rho)$ and $b\in(0,1)$, let $Q_b(\Tt)=OXVY$ denote the decorated quadrilateral with a Delzant corner at the origin $O$, edges $\overline{OX}$ and $\overline{OY}$ along the positive coordinate axes, nodal rays based at $X$ and $Y$ with directions
\[
\vec n_X=(p_\rho-6q_\rho,q_\rho),
\qquad
\vec n_Y=(q_\la,-p_\la),
\]
and axis lengths
\[
|OX|=3-b-\frac{d_\rho-bm_\rho}{q_\rho},
\qquad
|OY|=\frac{d_\la-bm_\la}{q_\la}.
\]
The remaining edge directions and lengths, and the nodal ray at $V$, are determined by $\Tt$ as in \cite[Lemma~2.1.4 and Lemma~2.3.3]{M2}. Following \cite[Definition~2.3.1]{M2}, we say that $Q_b(\Tt)$ is well-defined if its three nodal rays point into its interior, $|OX|,|OY|>0$, and $V$ lies in the interior of the positive quadrant. See Figure~\ref{fig:starting}.
\end{definition}

By \cite[Theorem 1.1.4]{M2}, the algebraic mutations agree
with the corresponding geometric mutations whenever
the latter are well-defined. In particular, if
$wQ_b^0$ is well-defined, then
\begin{equation}\label{eq:algGeom}
wQ_b^0=Q_b(w\Tt_0).
\end{equation}
The algebraic word $w\Tt_0$ can be defined even
when the geometric word $wQ_b^0$ is not well-defined
at the chosen value of $b$.

 \subsection{Background on staircases} \label{ss:stair}

 Here, we give the background needed for the staircase parameters with principal
and nonprincipal staircases. We review their description in terms of recursive triples from \cite{MMW,MPW,M2}.

By the classification in \cite{MMW,MPW}, for each staircase parameter $b\neq0,1/3$ there are a finite algebraic mutation word $w_b$ in $x,y,\ov{x},\ov{y},s,\ov{s}$ and an infinite sequence $a_{b,1},a_{b,2},\ldots\in\{x,y\}$ such that $w_b\Tt_0$ is one of the triples $S^iR^\de(\mathcal{B}_{2n})$ described below, and the triples obtained from $w_b\Tt_0$ by applying $a_{b,1},a_{b,2},\ldots$ in turn generate the staircase of $b$. (For the family starting from $R(\mathcal{B}_0)$, we modify this choice below.)

\begin{definition}\label{def:staircaseData}
Let $b\neq0,1/3$ be a staircase parameter, with $w_b$ and $(a_{b,k})_{k\geq1}$ as above. We call $w_b$ the \emph{starting word}, $\Tt_b:=w_b\Tt_0$ the \emph{starting triple}, and $(a_{b,k})_{k\geq1}$ the \emph{staircase sequence} of $b$. We set
\[
v_{b,0}=\mathrm{id},
\qquad
v_{b,k}=a_{b,k}\cdots a_{b,1},
\qquad
\Tt_{b,k}:=v_{b,k}w_b\Tt_0
=(\bE_{\la,k},\bE_{\mu,k},\bE_{\rho,k}).
\]
When $b$ is fixed, we suppress it from the notation and write $w$, $\Tt$, $a_k$, $v_k$, and $\Tt_k$.
\end{definition}

By \cite{MMW,MPW}, the middle entries satisfy
\[
\frac{m_{\mu,k}}{d_{\mu,k}}\longrightarrow b,
\qquad
\frac{p_{\mu,k}}{q_{\mu,k}}\longrightarrow\acc(b)
\]
where here $(d_{\mu,k},m_{\mu,k})$ are the degree coordinates defined in \eqref{eq:degree}. 

A staircase is called principal if its staircase sequence is eventually constant, that is, from some point on the subsequent
mutations are all $x$ or all $y$. Otherwise, it is
called nonprincipal; in this case the sequence contains infinitely many $x$-mutations and infinitely many $y$-mutations. By \cite[Section~2.2]{MPW}, this agrees with the description at the start of this section: if the sequence is eventually all $y$, then $c_b$ has only an ascending staircase; if it is eventually all $x$, then $c_b$ has only a descending staircase; and if it is nonprincipal, then $c_b$ has both.

In \cite[Section 3.3]{M2}, we gave intervals of parameters $b'$
for which the finite words $w_bQ_{b'}^0$ are
well-defined. For an ascending principal staircase, choose $N\geq0$ so that all mutations after the first $N$ are $y$, and set $w:=v_{b,N}w_b$, absorbing the finite initial segment into $w$. The results of \cite{M1,M2} show that $y^kwQ_b^0$ is well-defined for every $k\geq0$ and that these diagrams converge to a triangle giving a full filling at $\acc(b)$.

Here, we consider the descending principal
staircases and the nonprincipal staircases,
whose staircase sequences contain infinitely
many $x$-mutations. We will show that, with the starting words $w_b$ of Definition~\ref{def:staircaseData} (modified for the family starting from $R(\mathcal{B}_0)$ as described below),
these sequences also give diagrams converging
to a full-filling triangle at $\acc(b)$.
To do so, we first recall some properties of
the recursive triples and their associated
diagrams.

As described in \cite[Proposition 3.2.2]{M2}, the basic starting
triples are $\mathcal{B}_{2n}$, $n\geq0$, realized by
\[
\mathcal{B}_{2n}=y\ov{x}^{n+2}\Tt_0.
\]
The other starting families are obtained by applying
the shift and reflection symmetries $S$ and $R$.
Recall the map $S$ from \eqref{eq:S}, and set
\[
R(p,q,t)=(6p-35q,p-6q,-t).
\]
These symmetries act on triples by
\[
S(\Tt)=(S\bE_\la,S\bE_\mu,S\bE_\rho),
\qquad
R(\Tt)=(R\bE_\rho,R\bE_\mu,R\bE_\la).
\]
Both $S$ and $R$ change the sign $\eps$ to $-\eps$.
The reversal of the entries in $R(\Tt)$ accounts
for the reversal of the order of their centers.
Thus the starting families are
\[
S^iR^\de(\mathcal{B}_{2n}),
\qquad
i,n\geq0,\quad \de\in\{0,1\}.
\]
We refer to the families with $(i,\de)\neq(0,0)$
as the \emph{shifted} and \emph{reflected} families.
Finite mutation words realizing these families
are given in \cite[Proposition 3.2.2]{M2}. All other triples
in the staircase sequences are obtained by
applying finite sequences of $x$- and $y$-mutations
to these starting triples.

We recall the sufficient conditions from \cite[Section~3.3]{M2} for the corresponding geometric mutations to be well-defined. The finite words realizing the basic, shifted, and reflected starting triples are given in \cite[Proposition~3.2.2]{M2}, with the replacement for $R(\mathcal{B}_0)$ described in \cite[Lemma~3.2.3 and Corollary~3.2.4]{M2}. The conditions below apply to words of the form $w=u\,w_b$, where $w_b$ is such a starting word and $u$ is any finite word in $x,y$. For such a word $w$, write
\[
w\Tt_0=(\bE_\la,\bE_\mu,\bE_\rho),
\]
and let $\eps$ be the sign of this triple.
For the shifted and reflected families, including
the replacement for $R(\mathcal{B}_0)$, every mutation in
$wQ_{b'}^0$ is well-defined provided that
\[
\frac{\eps}{3}<\eps b'<
\eps\frac{m_\la}{d_\la};
\]
see \cite[Lemmas~3.3.13--3.3.16]{M2}.
For the basic family, where $w_b=y\ov{x}^{n+2}$, the word has the form
\[
w=u\,y\ov{x}^{n+2}.
\]
By \cite[Lemma~3.3.11]{M2}, every mutation in
$wQ_{b'}^0$ is well-defined provided that
\[
\frac{n+1}{n+2}<b'<\frac{m_\la}{d_\la}.
\]

To apply these results, we must show that the
staircase parameter $b$ lies in the corresponding
interval. For this, we use the following condition
on recursive triples, introduced in \cite[Section 4.2]{M2}.

\begin{definition}
A recursive triple
$\Tt=(\bE_\la,\bE_\mu,\bE_\rho)$ with sign $\eps$
is \textbf{correctly ordered} if
\[
3+2\sqrt{2}
<
\frac{p_\la}{q_\la}
<
\frac{p_\mu}{q_\mu}
<
\frac{p_\rho}{q_\rho},
\]
and, for $\alpha\in\{\la,\rho\}$,
\[
(0,0,2)
<
(p_\alpha,q_\alpha,\eps t_\alpha)
\leq
(p_\mu,q_\mu,\eps t_\mu),
\]
where the vector inequalities are interpreted
coordinatewise.
\end{definition}
Here $3+2\sqrt2=\acc(1/3)$, which is the minimum value of $\acc(b)$ for $b\in[0,1)$.

\begin{lemma}\label{lem:orderedSeeds}
The triples $S^iR^\de(\mathcal{B}_{2n})$,
where $i,n\geq0$ and $\de\in\{0,1\}$,
are correctly ordered, with the exception of
$R(\mathcal{B}_0)$. Moreover, $xy^rR(\mathcal{B}_0)$ is correctly ordered for every $r\geq0$.
\end{lemma}
\begin{proof}
We use the formulas for the basic triples in
\cite[Example~2.1.7]{MMW} and the symmetry formulas
in \cite[Section~2.3]{MMW}, expressed in our signed
$t$-coordinate notation.

The basic triples are
\[
\mathcal B_{2n}
=(\bE_\lambda,\bE_\mu,\bE_\rho),
\qquad n\geq0,
\]
where, in our signed $(p,q,t)$-coordinates,
\[
\begin{aligned}
\bE_\lambda&=(2n+6,\,1,\,2n+3),\\
\bE_\mu&=(4n^2+22n+29,\,2n+4,\,4n^2+16n+13),\\
\bE_\rho&=(2n+8,\,1,\,2n+5).
\end{aligned}
\]
These triples have sign $\eps=1$; see
\cite[Example~2.1.7]{MMW}. Applying $R$ gives
\[
\begin{aligned}
R(\mathcal B_{2n})=\bigl(&(12n+13,\,2n+2,\,-(2n+5)),\\
&(24n^2+62n+34,\,4n^2+10n+5,\,-(4n^2+16n+13)),\\
&(12n+1,\,2n,\,-(2n+3))\bigr),
\end{aligned}
\]
with sign $-1$. The definition can be checked directly for
$\mathcal B_{2n}$ and, when $n\geq1$, for
$R(\mathcal B_{2n})$. The triple $R(\mathcal B_0)$
has an entry $(1,0,-3)$ and is therefore not
correctly ordered.

It remains to check that the inequalities persist
under shifts. Set $c=3+2\sqrt{2}$. On $z=p/q$,
$S$ acts by
\[
z\longmapsto 6-\frac1z,
\]
which is strictly increasing and preserves
$(c,\infty)$. Hence shifts preserve the required
ordering of the ratios $p/q$.

For the coordinate comparisons, the
formulas above for $\mathcal B_{2n}$ give
\[
(p_\mu-p_\lambda,q_\mu-q_\lambda)
=(4n^2+20n+23,2n+3),
\]
and
\[
(p_\mu-p_\rho,q_\mu-q_\rho)
=(4n^2+20n+21,2n+3).
\]
For the reflected triple \(R(\mathcal B_{2n})\), these differences are transformed by \(R\), with the outer indices interchanged.
Each difference $(a,b)$ satisfies $a>6b>0$.
Under reflection, it becomes
\[
(6a-35b,a-6b),
\]
whose first coordinate is larger than its second,
and both are positive. Thus, for both the basic
and reflected triples,
\[
p_\mu-p_\alpha>q_\mu-q_\alpha>0,
\qquad \alpha\in\{\lambda,\rho\}.
\]
This inequality is preserved by $S$, since
\[
a>b>0 \quad\Longrightarrow\quad 6a-b>a>0.
\]
Hence the required coordinate comparisons hold
after every shift. Moreover, $S$ changes both
$\eps$ and $t$ by a sign, so the quantities
$\eps t_\alpha$ are unchanged.

Finally,
\[
R(\mathcal B_0)
=((13,2,-5),(34,5,-13),(1,0,-3)).
\] Applying $S$ to the exceptional entry
$(1,0,-3)$ gives $(6,1,3)$. The other entries of
$R(\mathcal{B}_0)$ have finite centers greater
than $6$, so their images have centers strictly
between $c$ and $6$. Consequently,
$SR(\mathcal{B}_0)$ is correctly ordered, and
the same argument applies to all its further
shifts.

It remains to show that $xy^rR(\mathcal B_0)$ is correctly ordered for every $r\geq0$. Since $R(\mathcal B_0)$ has sign $-1$ and right entry $(1,0,-3)$, each $y$-mutation acts by $(\bE_\la,\bE_\mu,(1,0,-3))\mapsto(\bE_\mu,3\bE_\mu-\bE_\la,(1,0,-3))$. Hence
\[
y^rR(\mathcal B_0)=(\bE_r,\bE_{r+1},(1,0,-3)),
\]
where $\bE_0=(13,2,-5)$, $\bE_1=(34,5,-13)$, and $\bE_{k+1}=3\bE_k-\bE_{k-1}$. Writing $\bE_k=(p_k,q_k,t_k)$, induction shows that $p_k$, $q_k$, and $-t_k$ are positive and increasing, that $-t_k>q_k$, and that $p_{k+1}q_k-p_kq_{k+1}=3$. The $x$-mutation then gives
\[
xy^rR(\mathcal B_0)=\bigl(\bE_r,\ -t_r\bE_{r+1}-(1,0,-3),\ \bE_{r+1}\bigr),
\]
with sign $-1$. Write $(p,q,t)$ for its middle entry. The ratios $p_k/q_k$ increase from $13/2>3+2\sqrt2$, and
\[
\frac pq=\frac{p_{r+1}}{q_{r+1}}-\frac1{-t_rq_{r+1}},
\qquad
\frac{p_{r+1}}{q_{r+1}}-\frac{p_r}{q_r}=\frac3{q_rq_{r+1}}>\frac1{-t_rq_{r+1}},
\]
using $-t_r>q_r$. Hence $p_r/q_r<p/q<p_{r+1}/q_{r+1}$. Since $-t_r\geq5$, we have $p=-t_rp_{r+1}-1>p_{r+1}$ and $q=-t_rq_{r+1}>q_{r+1}$, and $-t=(-t_r)(-t_{r+1})-3>-t_{r+1}>2$. Thus $xy^rR(\mathcal B_0)$ is correctly ordered. This proves the lemma.
\end{proof}

For the family starting from $R(\mathcal{B}_0)$, we modify the starting word, since $R(\mathcal{B}_0)$ has $q_\rho=0$ and so is not correctly ordered. If the sequence following $R(\mathcal{B}_0)$ contains an $x$, as for descending principal and nonprincipal staircases, let $r\geq0$ be the number of $y$-mutations before the first one. By \cite[Lemma~3.2.3]{M2},
\[
xy^rR(\mathcal{B}_0)=y\ov{x}^{r+2}\ov{y}s\,\Tt_0.
\]
In this case we set $w_b:=y\ov{x}^{r+2}\ov{y}s$ and $\Tt_b:=xy^rR(\mathcal{B}_0)$, and the staircase sequence $(a_{b,k})$ consists of the mutations following this first $x$. This does not affect the limits above, since the triples $\Tt_{b,k}$ form a tail of the original sequence.

\begin{cor}\label{cor:orderedStaircase}
Let $b$ be a staircase parameter whose staircase is descending principal or nonprincipal, with starting triple $\Tt_b=w_b\Tt_0$ and staircase sequence $(a_{b,k})$. Then $\Tt_b$ is correctly ordered, and so is $\Tt_{b,k}=v_{b,k}w_b\Tt_0$ for every $k\geq0$ where $w_b,v_{b,k}$ are defined in Definition~\ref{def:staircaseData}. 
\end{cor}
\begin{proof}
If $\Tt_b$ is one of the triples $S^iR^\de(\mathcal{B}_{2n})$, other than $R(\mathcal{B}_0)$, it is correctly ordered by Lemma~\ref{lem:orderedSeeds}. In the family starting from $R(\mathcal{B}_0)$, the staircase sequence contains an $x$, so $\Tt_b=xy^rR(\mathcal{B}_0)$ as above, which is correctly ordered by Lemma~\ref{lem:orderedSeeds}. Since $x$- and $y$-mutations preserve correct ordering \cite[Lemma~4.2.3]{M2}, every $\Tt_{b,k}$ is correctly ordered.
\end{proof}

 \subsection{Descending principal and nonprincipal staircases} \label{ss:descending}

In this subsection we treat the staircase parameters whose staircases are descending principal or nonprincipal. For each such $b$, we show in Lemma~\ref{lem:staircaseSeed} that every mutation in the explicit sequence of diagrams $Q_{b,k}$ defined below is well-defined, and in Lemma~\ref{lem:triangleLimit} that these diagrams converge nontrivially to a triangle. Before doing so, we establish two properties of the triples along the staircase sequence: Lemma~\ref{lem:diophLimits} computes the limits of their entries, and Lemma~\ref{lem:det} gives the inequality on their middle entries that we need for the mutations to be well-defined.

For a staircase parameter $b$, the diagrams we use are
\[
Q_{b,k}:=y\,v_{b,k}\,w_b\,Q_b^0,\qquad k\geq0,
\]
which we write as $Q_k=yv_kwQ_b^0$ when $b$ is fixed: we follow the staircase sequence and then apply one additional $y$-mutation. The reason is as follows. By \eqref{eq:algGeom} and Definition~\ref{def:Qb}, $|OY|=(d_\la-bm_\la)/q_\la$ depends only on the left entry of the triple, and an $x$-mutation does not change the left entry. In particular, along a descending principal staircase, where every subsequent mutation is an $x$, the diagrams $v_kwQ_b^0$ all have the same $|OY|$. The final $y$-mutation replaces the left entry with the middle entry $\bE_{\mu,k}$, whose ratio $m_{\mu,k}/d_{\mu,k}$ converges to $b$; this is what makes $|OY|$ converge to the correct value in Lemma~\ref{lem:triangleLimit}. It also means that the admissibility condition for $yv_kwQ_b^0$ is a condition on $\bE_{\mu,k}$, namely $\eps b<\eps\, m_{\mu,k}/d_{\mu,k}$, which is exactly what Lemma~\ref{lem:det} provides. The sequences $v_kwQ_b^0$, without the final $y$-mutation, were considered in \cite[Example~4.3.3]{M2}.

\begin{example}\label{ex:B0descending}
Consider the staircase parameter $b$ whose starting triple is the basic triple
\[
\Tt_b=\mathcal B_0=y\ov{x}^2\Tt_0=((6,1,3),(29,4,13),(8,1,5)),
\]
of sign $1$, and whose staircase sequence is $a_{b,k}=x$ for all $k$, so that $w_b=y\ov{x}^2$ and $v_{b,k}=x^k$. Since $\eps t_\la=3$ and the left entry is fixed, the middle entries satisfy \[\bE_{\mu,k+1}=3\bE_{\mu,k}-\bE_{\mu,k-1}\], giving
\[ \frac{29}{4}, \frac{79}{11}, \frac{208}{29}, \frac{545}{76}, \ldots ,\frac{p_{\mu,k}}{q_{\mu,k}}, \ldots, \to \acc(b) \approx 7.17.
\]
with degree coordinates $(d,m)$ as $(14,9),(38,24),(100,63),(262,165),\ldots$. The ratios $m_{\mu,k}/d_{\mu,k}$ converge to $b\approx0.6297$, which lies in the range $b>1/2$ used in the proof of Lemma~\ref{lem:staircaseSeed}, so the diagrams $Q_{b,k}$ are well-defined. 

The diagrams $Q_{b,k}=y\,x^k\,y\ov{x}^2\,Q_b^0$ have nodal ray at $Y$ of slope
\[
-\frac{29}{4},\ -\frac{79}{11},\ -\frac{208}{29},\ \ldots,
\]
recording the $p_k/q_k$ of the staircase. 
Without the final $y$-mutation, the left entry would remain $(6,1,3)$ and $|OY|$ would be constant. See Figure~\ref{fig:descending} for the upper bounds supplied by visible embeddings, approaching the full-filling point. 
\end{example}

\begin{figure}
    \centering
\begin{tikzpicture}[scale=0.75,x=1cm,y=1cm,>=Stealth,line cap=round,line join=round,
  every node/.style={font=\large}]
\draw[volume,line width=1.2pt,->] (0,0)--(10.2,0) node[right] {$z$};
\draw[volume,line width=1.2pt,->] (0,0)--(0,7.6) node[above] {$c_b$};
\draw[volume,line width=1.5pt] (0,0) .. controls (1,1.8) and (1.3,2.05) .. (1.9,2.5)
 .. controls (4.2,4.15) and (7.3,4.95) .. (9.7,5.25)
 node[right] {$V_b(z)$};
\draw[red!80!black,dashed] (1.9,0)--(1.9,2.5);
\fill[red!80!black] (1.9,2.5) circle (2.7pt);
\node[below=8pt,text=red!80!black] at (1.9,0) {$\operatorname{acc}(b)$};
\foreach \x/\y/\i in {3.15/4.1/2,5.1/5.5/1,8.05/7.15/0}{
  \draw[guide!65,densely dashed,line width=.55pt] (\x,0)--(\x,\y);
  \node[below=8pt,text=guide] at (\x,0) {$\dfrac{p_{\i}}{q_{\i}}$};
}
\draw[stairs,line width=1.5pt] (2.71,3.32)
 --(3.15,4.1)--(3.9,4.1)--(5.1,5.5)--(6.45,5.5)--(8.05,7.15)--(9.25,7.15);
\draw[embedding,dashed,line width=1.1pt] (3.15,4.1)--(3.282,4.334);
\draw[embedding,dashed,line width=1.1pt] (5.1,5.5)--(5.46,5.92);
\draw[embedding,dashed,line width=1.1pt] (8.05,7.15)--(9.17,8.305);
\foreach \p in {(3.15,4.1),(5.1,5.5),(8.05,7.15)}
  \fill[stairs] \p circle (3.2pt);
\foreach \p in {(3.282,4.334),(5.46,5.92),(9.17,8.305)}
  \fill[embedding] \p circle (2.8pt);
\node[embedding,above left=3pt] at (3.282,4.334) {$E(A_2,B_2)$};
\node[embedding,above left=3pt] at (5.46,5.92) {$E(A_1,B_1)$};
\node[embedding,above right=3pt] at (9.17,8.305) {$E(A_0,B_0)$};
\end{tikzpicture}
    \caption{A schematic of a descending staircase. We illustrate the case where the classes $\bE_k$ give outer corners of the staircase. Suppose $b$ has a descending staircase with starting triple $\Tt_b=w_b\Tt_0$. The triples
$\Tt_k:=x^kw_b\Tt_0=(\bE,\bE_k,\bE_{k-1})$
determine outer corners at $z=p_k/q_k$, with $p_k/q_k\to\acc(b)$ as $k\to\infty$.
The diagrams $Q_k:=yx^kw_bQ_b^0=OX_kV_kY_k$
supply upper bounds via visible embeddings, labeled $E(A_k,B_k)$, where
$A_k=|OX_k|$ and $B_k=|OY_k|$.
The point labeled $E(A_k,B_k)$ has coordinates $\left(\max\{A_k/B_k,B_k/A_k\},\;1/\min\{A_k,B_k\}\right)$.
The corresponding purple points lie strictly above the graph
of the ellipsoid embedding function, on the extensions of
the obstruction lines associated with $\bE_k$.
As $k\to\infty$, the distance between each purple point and
its corresponding orange outer corner tends to zero,
and these upper bounds approach the full-filling point at
the accumulation point.}
    \label{fig:descending}
\end{figure}

\begin{rmk} \label{rmk:descendVisible}
    For an ascending principal staircase, absorb the finite initial segment into a word $w$ so that the remaining sequence consists only of $y$-mutations, as in Section~\ref{ss:stair}. In \cite[Theorem~D and Theorem~1.3.3]{M2}, we showed that the diagrams $xy^kwQ_b^0$ give visible obstructions for all sufficiently large $k$, in the sense of \cite[Section~5.1]{McDuffSiegel2025}. Thus $y^kwQ_b^0$ gives the triangular limit, while the additional $x$-mutation is used to obtain visible obstructions. Hence, these diagrams can be used to show that $b$ has an infinite staircase. We discussed in \cite[Example~4.3.3]{M2} why the same method does not work for the descending principal staircases. We still do not know if visible obstructions can be used to construct the descending staircases using other sequences. 
\end{rmk}

 \begin{lemma}\label{lem:diophLimits}
Let $(d_k,m_k;p_k,q_k)$ be integers satisfying
\[
3d_k-m_k=p_k+q_k,\qquad
d_k^2-m_k^2=p_kq_k-1,
\]
with $d_k>0$ and $p_k>q_k>0$.
Assume $m_k/d_k\to b\in(0,1)$ and $q_k\to\infty$.
Then, setting $z_\infty=\acc(b)$ and
$D=(z_\infty+1)/(3-b)$, we have
\[
\frac{p_k}{q_k}\to z_\infty,\qquad
\frac{d_k}{q_k}\to D,\qquad
\frac{d_k-bm_k}{q_k}\to\frac{z_\infty}{D},
\]
where $(1-b^2)D^2=z_\infty$.
\end{lemma}

\begin{proof}
Set $b_k=m_k/d_k$ and $z_k=p_k/q_k$. We have
\[
3-b_k=\frac{3d_k-m_k}{d_k},
\qquad
1-b_k^2=\frac{d_k^2-m_k^2}{d_k^2}.
\]
Consequently,
\[
\frac{(3-b_k)^2}{1-b_k^2}
=
\frac{(3d_k-m_k)^2}{d_k^2-m_k^2}.
\]
Substituting the Diophantine equations
$3d_k-m_k=p_k+q_k$ and
$d_k^2-m_k^2=p_kq_k-1$ gives
\[
\frac{(3-b_k)^2}{1-b_k^2}
=\frac{(p_k+q_k)^2}{p_kq_k-1}.
\]
Therefore
\[
\frac{(3-b_k)^2}{1-b_k^2}
\left(1-\frac1{p_kq_k}\right)
=\frac{(p_k+q_k)^2}{p_kq_k}
=z_k+\frac1{z_k}+2.
\]
If $b_k\to b$ and $q_k\to\infty$, taking limits yields
\[
z_k+\frac1{z_k}
\longrightarrow
\frac{(3-b)^2}{1-b^2}-2.
\]
Since $z_k>1$, we obtain
$z_k\to\acc(b)$.

Rearranging the linear Diophantine equation gives
\[
\frac{d_k}{q_k}
=\frac{p_k/q_k+1}{3-m_k/d_k}
\longrightarrow
\frac{\acc(b)+1}{3-b}=D.
\]
Consequently,
\[
\frac{d_k-bm_k}{q_k}
=\frac{d_k}{q_k}\left(1-b\frac{m_k}{d_k}\right)
\longrightarrow D(1-b^2).
\]
Dividing the quadratic Diophantine equation by $q_k^2$
and taking limits gives
\[
D^2(1-b^2)=z_\infty.
\]
Therefore $D(1-b^2)=z_\infty/D$, proving the final limit.
\end{proof}

\begin{lemma}\label{lem:det}
Let $\Tt$ be a correctly ordered recursive triple
of sign $\eps$, and let $\bE_{\mu,k}$ be the middle
entry after $k$ mutations in a sequence of $x$-
and $y$-mutations starting from $\Tt$.
Then $\eps m_{\mu,k}/d_{\mu,k}$ is strictly decreasing.
If $m_{\mu,k}/d_{\mu,k}\to b$, then
\[
\eps(m_{\mu,k}-bd_{\mu,k})>0
\qquad\text{for all }k\geq0.
\]
\end{lemma}

\begin{proof}
Let $a_1,a_2,\ldots\in\{x,y\}$ be the mutation sequence.
Set $v_0=\mathrm{id}$ and $v_k=a_k\cdots a_1$, and write
\[
v_k\Tt=(\bE_{\la,k},\bE_{\mu,k},\bE_{\rho,k}).
\]
By \cite[Lemma~4.2.3]{M2}, each triple remains correctly ordered and
has sign $\eps$. Moreover, the determinant identities
in \cite[Lemma~3.3.5]{M2} give
\[
m_{\la,k}d_{\mu,k}-m_{\mu,k}d_{\la,k}
=\eps q_{\rho,k},
\qquad
m_{\rho,k}d_{\mu,k}-m_{\mu,k}d_{\rho,k}
=\eps p_{\la,k}.
\]

If $a_{k+1}=x$, then
\[
\bE_{\mu,k+1}
=\eps t_{\la,k}\bE_{\mu,k}-\bE_{\rho,k}.
\]
Since the degree coordinates satisfy the same
linear relation, we obtain
\[
m_{\mu,k}d_{\mu,k+1}-m_{\mu,k+1}d_{\mu,k}
=m_{\rho,k}d_{\mu,k}-m_{\mu,k}d_{\rho,k}
=\eps p_{\la,k}.
\]
Similarly, if $a_{k+1}=y$, then
\[
\bE_{\mu,k+1}
=\eps t_{\rho,k}\bE_{\mu,k}-\bE_{\la,k},
\]
so
\[
m_{\mu,k}d_{\mu,k+1}-m_{\mu,k+1}d_{\mu,k}
=m_{\la,k}d_{\mu,k}-m_{\mu,k}d_{\la,k}
=\eps q_{\rho,k}.
\]
Thus
\[
\eps\left(
\frac{m_{\mu,k}}{d_{\mu,k}}
-\frac{m_{\mu,k+1}}{d_{\mu,k+1}}
\right)
=
\begin{cases}
\displaystyle
\frac{p_{\la,k}}{d_{\mu,k}d_{\mu,k+1}}>0,
&a_{k+1}=x,\\[8pt]
\displaystyle
\frac{q_{\rho,k}}{d_{\mu,k}d_{\mu,k+1}}>0,
&a_{k+1}=y.
\end{cases}
\]
Hence $\eps m_{\mu,k}/d_{\mu,k}$ is strictly decreasing.
If it converges to $\eps b$, every term is strictly
larger than $\eps b$. Since $d_{\mu,k}>0$, this gives
\[
\eps(m_{\mu,k}-bd_{\mu,k})>0
\]
for every $k\geq0$.
\end{proof}

We now show that the limiting $b$-values corresponding
to descending principal and nonprincipal staircases
satisfy the sufficient conditions recalled in Section~\ref{ss:stair}.

\begin{lemma}\label{lem:staircaseSeed}
Let $b\in(0,1)$ be a staircase parameter whose staircase is descending principal or nonprincipal. Then the diagrams
\[
Q_{b,k}=y\,v_{b,k}\,w_b\,Q_b^0
\]
are well-defined for every $k\geq0$ where $v_{b,k},w_b$ are defined as in Definition~\ref{def:staircaseData}. 
\end{lemma}

\begin{proof}
Fix $b$ and suppress it from the notation, so that $\Tt=w_b\Tt_0$ and $\Tt_k=v_{b,k}w_b\Tt_0$ as in Definition~\ref{def:staircaseData}. By Corollary~\ref{cor:orderedStaircase}, $\Tt$ and every $\Tt_k$ are correctly ordered.

Let $\eps$ be the sign of $\Tt$. By Section~\ref{ss:stair},
\[
\frac{m_{\mu,k}}{d_{\mu,k}}\longrightarrow b,
\]
and we use this limit to check the sufficient conditions
for every finite word $yv_kw$ to be well-defined.

First, Lemma~\ref{lem:det} gives
\[
\eps b<\eps\frac{m_{\mu,k}}{d_{\mu,k}}
\qquad\text{for every }k\geq0.
\]
Next, since each triple remains correctly ordered,
$\eps t_{\mu,k}>0$. Using
$t_{\mu,k}=3m_{\mu,k}-d_{\mu,k}$, we obtain
\[
\eps\left(
3\frac{m_{\mu,k}}{d_{\mu,k}}-1
\right)
=\frac{\eps t_{\mu,k}}{d_{\mu,k}}>0.
\]
Taking limits gives $\eps(3b-1)\geq0$.
Since $b\neq1/3$ for the staircases under consideration,
this inequality is strict. Thus
\begin{equation}\label{eq:seedInterval}
\frac{\eps}{3}<\eps b<
\eps\frac{m_{\mu,k}}{d_{\mu,k}}
\qquad\text{for every }k\geq0.
\end{equation}

The final $y$-mutation gives
\[
yv_k\Tt
=(\bE_{\mu,k},\ \eps t_{\rho,k}\bE_{\mu,k}-\bE_{\la,k},\ \bE_{\rho,k}),
\]
so its left entry is $\bE_{\mu,k}$.
For the shifted and reflected families, including the family starting from $R(\mathcal{B}_0)$, the word $yv_kw$ has the form $u\,w_b$ with $u=yv_k$, so by \cite[Lemmas~3.3.13--3.3.16]{M2}, as recalled in Section~\ref{ss:stair}, the required admissibility condition is
\[
\frac{\eps}{3}<\eps b<
\eps\frac{m_{\mu,k}}{d_{\mu,k}},
\]
which holds by \eqref{eq:seedInterval}.
Hence every mutation in $yv_kwQ_b^0$ is well-defined.

For the family generated by $\mathcal{B}_{2n}$, the starting word is $w=y\ov{x}^{n+2}$, so
\[
yv_kw=(yv_k)\,y\ov{x}^{n+2}
\]
has the form $u\,y\ov{x}^{n+2}$ with $u=yv_k$ a finite word in $x,y$, and the admissibility criterion
in \cite[Lemma 3.3.11]{M2} requires
\[
\frac{n+1}{n+2}<b<
\frac{m_{\mu,k}}{d_{\mu,k}}.
\]
Here $\eps=1$, so the upper bound follows from
\eqref{eq:seedInterval}.
The classification places every staircase parameter
in this family above $(n+1)/(n+2)$; see \cite[Remark~2.1.8]{MMW}.
Hence the lower bound also holds, and every
$yv_kwQ_b^0$ is well-defined in this case as well.
\end{proof}

The following lemma adapts the triangular-limit
arguments in \cite{M1,M2}, which treated sequences
$v_k=y^k$, to sequences containing infinitely many
$x$-mutations.

\begin{lemma}\label{lem:triangleLimit}
Let $w$ be a finite mutation word such that $\Tt=w\Tt_0$ is correctly ordered, and let $a_1,a_2,\ldots\in\{x,y\}$ contain infinitely many $x$'s. With $v_k=a_k\cdots a_1$ and $\Tt_k=v_k\Tt=(\bE_{\la,k},\bE_{\mu,k},\bE_{\rho,k})$, assume that $m_{\mu,k}/d_{\mu,k}\to b\in(0,1)$ and that $Q_k=yv_kwQ_b^0$ is well-defined for every $k\geq0$. Then the standardized polygons of the $Q_k$ are pairwise distinct and converge nontrivially, in the sense of Definition~\ref{def:converge}, to a triangle giving a full filling at $\acc(b)$.
\end{lemma}

Lemma~\ref{lem:triangleLimit} does not use that $b$ is a staircase parameter or that $(a_k)$ is a staircase sequence; it applies to any such sequence of mutations whose middle entries have $m_{\mu,k}/d_{\mu,k}\to b$.

\begin{proof}
By \cite[Lemma~4.2.3]{M2}, each $\Tt_k$ remains correctly ordered.
Correct ordering and the mutation formulas give
\[
q_{\mu,k+1}\geq2q_{\mu,k},
\]
whether $a_{k+1}=x$ or $a_{k+1}=y$.
Indeed, if $a_{k+1}=x$ then $q_{\mu,k+1}=\eps t_{\la,k}q_{\mu,k}-q_{\rho,k}$, and correct ordering gives $\eps t_{\la,k}\geq3$ and $q_{\rho,k}\leq q_{\mu,k}$; the case $a_{k+1}=y$ is the same with the roles of $\la$ and $\rho$ exchanged.
Thus $q_{\mu,k}\to\infty$.

Since
\[
\frac{m_{\mu,k}}{d_{\mu,k}}\to b
\quad\text{and}\quad q_{\mu,k}\to\infty,
\]
Lemma~\ref{lem:diophLimits} gives
\[
\frac{p_{\mu,k}}{q_{\mu,k}}
\longrightarrow z_\infty:=\acc(b)>1.
\]

Each $x$-mutation replaces the right entry by the
previous middle entry, while each $y$-mutation preserves
the right entry. Since $x$ occurs infinitely often,
the right entries are middle entries from stages tending
to infinity. Consequently,
\[
\frac{m_{\rho,k}}{d_{\rho,k}}\to b,
\qquad q_{\rho,k}\to\infty.
\]

Applying Lemma~\ref{lem:diophLimits} to both the middle
and right entries, we obtain
\[
\frac{d_{\mu,k}-bm_{\mu,k}}{q_{\mu,k}}
\longrightarrow\frac{z_\infty}{D},
\qquad
\frac{d_{\rho,k}-bm_{\rho,k}}{q_{\rho,k}}
\longrightarrow\frac{z_\infty}{D},
\]
where
\[
D=\frac{z_\infty+1}{3-b}>0,
\qquad
(1-b^2)D^2=z_\infty.
\]

Write $Q_k=OX_kV_kY_k$, where
$X_k=(A_k,0)$ and $Y_k=(0,B_k)$.
By \eqref{eq:algGeom}, $Q_k=Q_b(yv_kw\Tt_0)=Q_b(y\Tt_k)$. Since the left and right entries of $y\Tt_k$ are
$\bE_{\mu,k}$ and $\bE_{\rho,k}$, respectively,
Definition~\ref{def:Qb} gives
\[
A_k=3-b-\frac{d_{\rho,k}-bm_{\rho,k}}{q_{\rho,k}},
\qquad
B_k=\frac{d_{\mu,k}-bm_{\mu,k}}{q_{\mu,k}}.
\]
Hence,
\[
B:=\lim_{k\to\infty}B_k=\frac{z_\infty}{D},
\qquad
A:=\lim_{k\to\infty}A_k
=3-b-\frac{z_\infty}{D}=\frac1D.
\]
It follows that
\[
A+B=3-b,\qquad AB=1-b^2,\qquad
\frac BA=z_\infty=\acc(b).
\]

To see that $Q_k$ converges to $T:=\operatorname{conv}\{O,(A,0),(0,B)\}$, write $V_k=(u_k,v_k)$ with $u_k,v_k\geq0$. Since mutations preserve area, $A_kv_k+B_ku_k=1-b^2$, and since $A_k\to A>0$ and $B_k\to B>0$, the points $V_k$ are bounded. Any limit of the $Q_k$ is then a convex set containing $T$ with area $(1-b^2)/2=AB/2$, the same as $T$, so it must be $T$ itself.

Finally, we compare axis lengths. An $x$-mutation shears the part of the diagram on the far side of the nodal line through $X$, which contains $V$ but not $O$. After the mutation, the image of the edge $XV$ is collinear with $OX$, so $X$ is no longer a vertex and the horizontal axis edge now extends to the image of $V$; hence $|OX|$ strictly increases. A $y$-mutation shears only the part on the far side of the nodal line through $Y$, which contains $V$ but not the edge $OX$, so it preserves $|OX|$, and by the same reasoning it strictly increases $|OY|$.
The additional final $y$-mutation also preserves
this length. Thus $A_k$ is nondecreasing, and
$A_{k+1}>A_k$ whenever $a_{k+1}=x$.
If $a_{k+1}=y$, then $Q_{k+1}=yQ_k$, so
$A_{k+1}=A_k$ and $B_{k+1}>B_k$.
Between any two indices, either an $x$ occurs,
giving a strict increase in $A_k$, or only $y$'s
occur, giving a strict increase in $B_k$.
Hence the standardized polygons are pairwise distinct.
The convergence is therefore nontrivial in the sense of Definition~\ref{def:converge}, and Theorem~\ref{thm:general} gives a full filling at $\acc(b)$.
\end{proof}

\begin{prop}\label{prop:descending}
Let $b\in(0,1)$ be a staircase parameter whose staircase is descending principal or nonprincipal. Then there is a sequence $Q_{b,k}$ of almost toric base diagrams for $H_b$ in standard position that converges nontrivially to a triangle.
\end{prop}

\begin{proof}
Consider the diagrams
\[
Q_{b,k}:=y\,v_{b,k}\,w_b\,Q_b^0,\qquad k\geq0.
\]
These are well-defined by Lemma~\ref{lem:staircaseSeed}. Moreover, the starting triple $\Tt_b=w_b\Tt_0$ is correctly ordered by Corollary~\ref{cor:orderedStaircase}, the staircase sequence contains infinitely many $x$'s, and $m_{\mu,k}/d_{\mu,k}\to b$ (Section~\ref{ss:stair}). Hence Lemma~\ref{lem:triangleLimit}, applied with $w=w_b$ and $v_k=v_{b,k}$, shows that the $Q_{b,k}$ converge nontrivially to a triangle.
\end{proof}

 \subsection{Special rational $b$-values} \label{ss:specialRational}
\begin{figure}
    \centering
   \begin{tikzpicture}[scale=0.9,x=1cm,y=1cm,>=Stealth,line cap=round,
 every node/.style={font=\large}]
\draw[volume,line width=1.2pt,->] (0,0)--(10.5,0) node[right] {$z$};
\draw[volume,line width=1.2pt,->] (0,0)--(0,6.3) node[above] {$c_b$};
\draw[volume,line width=1.5pt] (0,0)
 .. controls (1.6,2.6) and (4.6,4.05) .. (8.3,4.8)
 .. controls (9.1,4.9622) and (9.8,5.08) .. (10.4,5.17)
 node[right] {$V_b(z)$};
\draw[embedding,dash pattern=on 0pt off 3.5pt,line width=1.15pt] (0,4.8)--(8.3,4.8);
\node[left=5pt] at (0,4.8) {$\frac52$};
\draw[red!80!black,dashed,line width=.7pt] (8.3,0)--(8.3,4.8);
\fill[red!80!black] (8.3,4.8) circle (2.8pt);
\node[below=7pt,text=red!80!black] at (8.3,0) {$\operatorname{acc}(b)=6$};
\foreach \x in {2.1,4.6,6.75}
 \fill[embedding] (\x,4.8) circle (2.8pt);
\node[embedding,above=10pt] at (2.1,4.8) {$E(A_0,B_0)$};
\node[embedding,above=10pt] at (4.6,4.8) {$E(A_1,B_1)$};
\node[embedding,above=10pt] at (6.75,4.8) {$E(A_2,B_2)$};
\end{tikzpicture}
    \caption{For $b=1/5$, the points labeled $E(A_k,B_k)$ are upper bounds
supplied by visible embeddings. The point labeled $E(A_k,B_k)$ has coordinates $\left(\max\{A_k/B_k,B_k/A_k\},\;1/\min\{A_k,B_k\}\right)$.
They lie on the horizontal line at height $5/2$, approaching
the full-filling point at $\acc(b)=6$ from the left. This is in contrast to the descending case of Figure~\ref{fig:descending}. This is not to scale.}
    \label{fig:b15}
\end{figure}
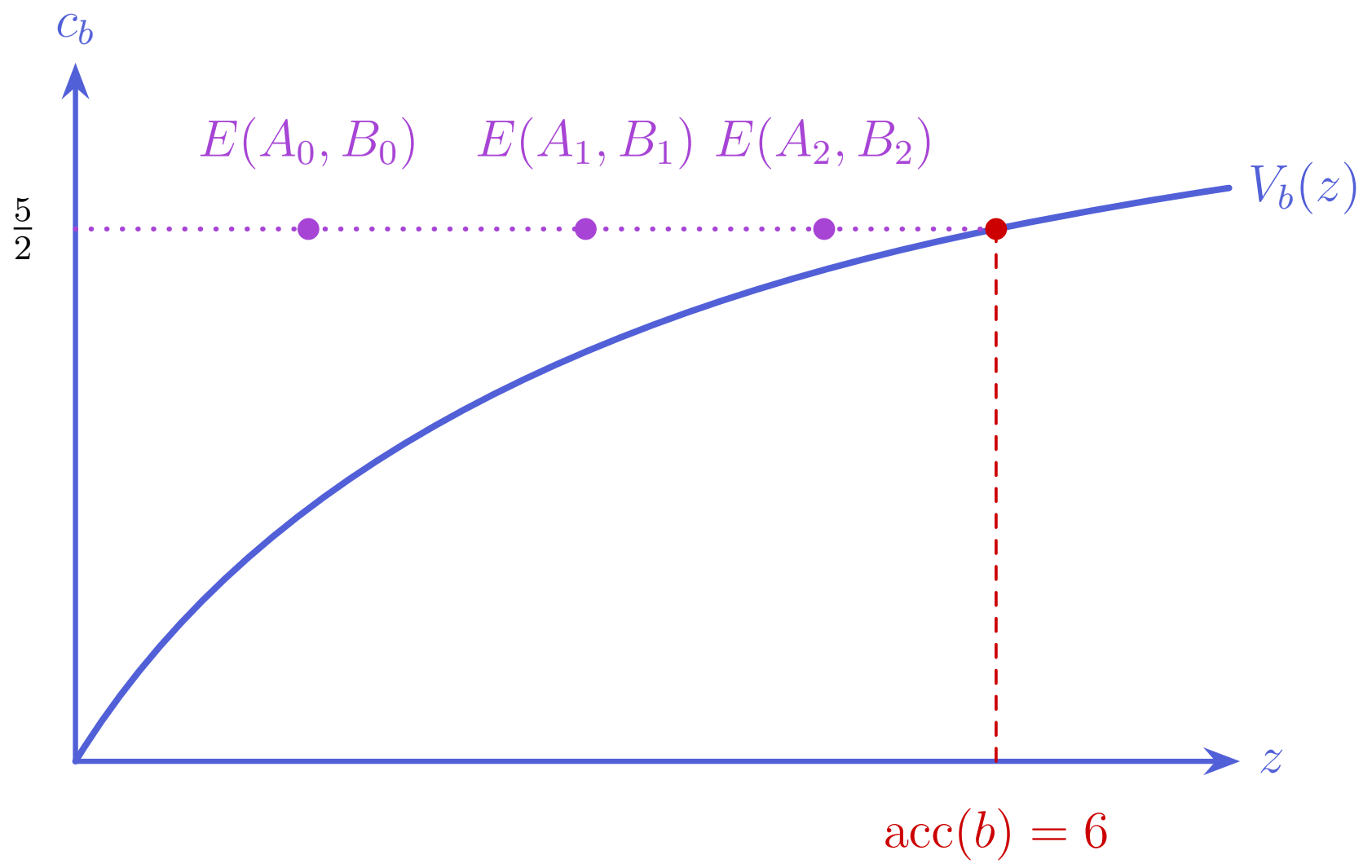

We now construct the required sequences for the special
rational values $b_i$, see \eqref{eq:special}. For each $i\geq0$,
we consider the diagrams
\[
Q_n=x^n\ov{s}^{3(i+1)}Q_{b_i}^0,
\qquad n\geq0.
\]
We will show in Proposition~\ref{lem:rationalTriangle} that every mutation in these words is
well-defined and that the resulting diagrams converge
nontrivially to a triangle. We discuss the case of $b_0=1/5$ in Example~\ref{ex:perimeterJump} and see Figure~\ref{fig:b15} for the upper bounds supplied by visible embeddings.

To describe the initial $\ov{s}$-mutations, define
\[
\bF_{-2}=(4,1,1),\qquad
\bF_{-1}=(2,1,-1),\qquad
\bF_0=(1,1,2),
\]
and
\[
\bF_{j+3}=S^{-1}\bF_j,
\qquad j\geq-2.
\]
The formula for $\ov{s}$ gives
\[
\ov{s}^{j}\Tt_0
=(\bF_j,\bF_{j-1},\bF_{j-2}),
\qquad j\geq0,
\]
with sign $(-1)^j$.
Write $\bF_j=(p_j,q_j,t_j)$, and let
$(d_j,m_j)$ be its degree coordinates.
Whenever $d_j\ne0$, set
\[
\alpha_j=\frac{m_j}{d_j}.
\]

We now give some useful numerical properties about the classes $\bF_j$ and their degree coordinates. 

\begin{lemma}\label{lem:rationalTriples}
For $i\geq0$, write $\mathfrak{p}=\mathfrak{p}_i$, $\mathfrak{q}=\mathfrak{q}_i$, and
$\eps=(-1)^{i+1}$. Then
\[
\begin{aligned}
\bF_{3i+3}&=(\mathfrak{p}-5\mathfrak{q},\mathfrak{p}-\mathfrak{q},2\eps),\\
\bF_{3i+2}&=(2\mathfrak{p}-11\mathfrak{q},\mathfrak{p}-2\mathfrak{q},-\eps),\\
\bF_{3i+1}&=(4\mathfrak{p}-23\mathfrak{q},\mathfrak{p}-4\mathfrak{q},\eps).
\end{aligned}
\]
\end{lemma}
\begin{proof}
Recall that
\[
S^{-1}(p,q,t)=(q,6q-p,-t)
\quad\text{and}\quad
\ov{s}^{3}\Tt=S^{-1}\Tt,
\]
where $S^{-1}$ acts entrywise. The recursion defining
$(\mathfrak{p}_i,\mathfrak{q}_i)$ gives
\[
S^{-(i+1)}(p,q,t)
=\bigl((\mathfrak{p}-6\mathfrak{q})p+\mathfrak{q}q,\,-\mathfrak{q}p+\mathfrak{p}q,\,\eps t\bigr).
\]
Applying this to
\[
\Tt_0=((1,1,2),(2,1,-1),(4,1,1))
\]
gives the claimed formulas. 
\end{proof}

The following lemma establishes positivity of $d_j$
and the nesting properties of these ratios.

\begin{lemma}\label{lem:alphaNesting}
For every $j\geq-2$, we have
\[
d_j>0,\qquad 0<p_j<6q_j.
\]
Moreover, the closed intervals with endpoints
$\alpha_j$ and $\alpha_{j+1}$ are nested as
$j\geq0$ increases:
\[
\alpha_1\leq \cdots \leq \alpha_{2j+1} \leq \cdots
<\frac13<
\cdots\leq\alpha_{2j}\leq\cdots\leq\alpha_0.
\]
\end{lemma}
\begin{proof}
Recall that
\[
(p_{j+3},q_{j+3},t_{j+3})
=(q_j,6q_j-p_j,-t_j).
\]
The initial values give
\[
(p_0,q_0)=(1,1),\qquad
(p_1,q_1)=(1,2),\qquad
(p_2,q_2)=(1,4).
\]
Since $(p,q)\mapsto(q,6q-p)$ preserves
$0<p\leq q$, these inequalities hold for every
$j\geq0$. Together with the initial values at
$j=-2,-1$, this proves $0<p_j<6q_j$ for all
$j\geq-2$. Moreover, $p_j,q_j$ are positive integers
and $|t_j|\leq2$, so
\[
d_j=\frac{3(p_j+q_j)+t_j}{8}>0.
\]

For the nesting assertions, write
\[
s^j\Tt_0
=(\mathbf{G}_j,\mathbf{G}_{j+1},\mathbf{G}_{j+2}),
\]
and let $J(p,q,t)=(q,p,t)$. The initial values satisfy
\[
\bF_0=J\mathbf{G}_0,\qquad
\bF_1=J\mathbf{G}_1,\qquad
\bF_2=J\mathbf{G}_2.
\]
Since $\mathbf{G}_{j+3}=S\mathbf{G}_j$ and
$S^{-1}J=JS$, the recursions imply
$\bF_j=J\mathbf{G}_j$ for every $j\geq0$.
Interchanging $p$ and $q$ preserves the degree
coordinates $(d,m)$, so our ratios $\alpha_j$
coincide with those in \cite[Lemma~3.3.9]{M2}.
That lemma gives
\[
(-1)^j\alpha_{j+1}
<\frac{(-1)^j}{3}
<(-1)^j\alpha_{j+2}
\leq(-1)^j\alpha_j.
\]
Taking $j$ even and odd gives the claimed
inequalities. Thus the lower endpoints of the
intervals are nondecreasing and the upper
endpoints are nonincreasing, proving the nesting.
\end{proof}

\begin{lemma}\label{lem:barsCriterion}
Let $\Tt=(\bE_\la,\bE_\mu,\bE_\rho)$ be a recursive triple of sign $\eps$, and assume that
$Q_{b}(\Tt)$ is well-defined.
Suppose that $p_\rho<6q_\rho$, and write
\[
S^{-1}\bE_\mu=(p,q,t)
\]
with degree coordinates $(d,m)$. Then:
\begin{enumerate}
\item $\ov{s}Q_b(\Tt)$ is well-defined if and only if $\eps(m-bd)>0$;
\item $xQ_b(\Tt)$ is well-defined if and only if $\eps(m-bd)<0$.
\end{enumerate}
In particular, if $d>0$, condition (i) is equivalent to
\[
\eps b<\eps\frac{m}{d}.
\]
\end{lemma}

\begin{proof}

Write $Q_b(\Tt)=OXVY$, and set $A=|OX|$ and $B=|OY|$.
A direct computation from the definitions of $S^{-1}$ and the degree coordinates gives
\[
d=3q_\mu-d_\mu,\qquad m=q_\mu-m_\mu.
\]
The numerator identity in the affine-length calculation
for the $x$-mutation in \cite[Lemma~2.3.10]{M2} gives
\begin{equation}\label{eq:barsIntersection}
\begin{aligned}
B(6q_\rho-p_\rho)-Aq_\rho
&=-\frac{\eps\bigl((m_\mu-q_\mu)-b(d_\mu-3q_\mu)\bigr)}{q_\la}\\
&=\frac{\eps(m-bd)}{q_\la}.
\end{aligned}
\end{equation}
The nodal ray from $X$ has direction
$(p_\rho-6q_\rho,q_\rho)$. Since $p_\rho<6q_\rho$,
it meets the positive $y$-axis at height
\[
h=\frac{Aq_\rho}{6q_\rho-p_\rho}.
\]
The $\ov{s}$-mutation is well-defined precisely when
this intersection lies below $Y$, that is, when $h<B$.
Since the nodal ray at $X$ leaves $Q_b(\Tt)$ through either $\overline{OY}$ or $\overline{VY}$, the $x$-mutation is well-defined precisely when $h>B$. These conditions are equivalent to
\[
B(6q_\rho-p_\rho)-Aq_\rho>0
\qquad\text{and}\qquad
B(6q_\rho-p_\rho)-Aq_\rho<0,
\]
respectively. By \eqref{eq:barsIntersection} and $q_\la>0$, they are equivalent to $\eps(m-bd)>0$ and $\eps(m-bd)<0$, respectively,
as desired.
\end{proof}

Similarly, the nodal ray at $Y$ leaves the diagram through either $\overline{OX}$ or $\overline{XV}$, and \cite[Lemma~3.3.1]{M2} gives the analogous criterion for the $s$- and $y$-mutations.

\begin{lemma}\label{lem:rationalSeed}
Let $b_i$ be a special rational value, and write
$\mathfrak{p}=\mathfrak{p}_i$, $\mathfrak{q}=\mathfrak{q}_i$, and $\eps=(-1)^{i+1}$.
Then every mutation in
\[
\ov{s}^{3(i+1)}Q_{b_i}^0
\]
is well-defined.
\end{lemma}

    \begin{proof}
Assume that $\ov{s}^{j}Q_b^0$ is well-defined.
Its associated triple is
$(\bF_j,\bF_{j-1},\bF_{j-2})$, with sign $(-1)^j$.
By Lemma~\ref{lem:alphaNesting},
\[
p_{j-2}<6q_{j-2},
\qquad
d_{j+2}>0.
\]
Since $S^{-1}\bF_{j-1}=\bF_{j+2}$,
Lemma~\ref{lem:barsCriterion} implies that the next
$\ov{s}$-mutation is well-defined if
\[
(-1)^j(\alpha_{j+2}-b)>0.
\]
The initial diagram $Q_b^0$ is well-defined for
$b\in(0,1)$, and $\alpha_0=1$, $\alpha_1=0$.
Since the intervals bounded by $\alpha_j$ and
$\alpha_{j+1}$ are nested by
Lemma~\ref{lem:alphaNesting}, induction shows that
$\ov{s}^{j}Q_b^0$ is well-defined whenever
\[
(-1)^j\alpha_{j+1}
<
(-1)^j b
<
(-1)^j\alpha_j.
\]
Thus, to prove the claim, it suffices to show that
\[
\eps\alpha_{3(i+1)+1}
<
\eps b_i
<
\eps\alpha_{3(i+1)},
\]
where $(-1)^{3(i+1)}=\eps$.

For the first inequality, $t_{3(i+1)+1}$ has sign
$-\eps$, whereas the definition of $b_i$ gives
\[
\eps\left(b_i-\frac13\right)
=\frac{8}{3\bigl(3(\mathfrak{p}+\mathfrak{q})+\eps\bigr)}>0.
\]
Hence
\[
\eps\alpha_{3(i+1)+1}<\frac{\eps}{3}<\eps b_i.
\]

For the second inequality, the formula from Lemma~\ref{lem:rationalTriples} 
\[
\bF_{3(i+1)}=(\mathfrak{p}-5\mathfrak{q},\mathfrak{p}-\mathfrak{q},2\eps)
\]
gives
\[
\alpha_{3(i+1)}
=\frac{\mathfrak{p}-3\mathfrak{q}+3\eps}{3(\mathfrak{p}-3\mathfrak{q})+\eps}.
\]
Consequently,
\[
\eps\left(\alpha_{3(i+1)}-\frac13\right)
=\frac{8}{3\bigl(3(\mathfrak{p}-3\mathfrak{q})+\eps\bigr)}
>
\frac{8}{3\bigl(3(\mathfrak{p}+\mathfrak{q})+\eps\bigr)}
=\eps\left(b_i-\frac13\right).
\]
Thus $\eps b_i<\eps\alpha_{3(i+1)}$, completing
the proof.
\end{proof}

\begin{prop}\label{lem:rationalTriangle}
Let $b_i$ be a special rational value. Then the diagrams
\[
Q_n=x^n\ov{s}^{3(i+1)}Q_{b_i}^0,
\qquad n\geq0,
\]
are well-defined, their standardized polygons are
pairwise distinct, and they converge nontrivially, in the sense of Definition~\ref{def:converge}, to a triangle
giving a full filling at $\acc(b_i)$.
\end{prop}
\begin{proof}
Fix $i$, and write
\[
\mathfrak{p}=\mathfrak{p}_i,\qquad \mathfrak{q}=\mathfrak{q}_i,\qquad
\eps=(-1)^{i+1},\qquad N=3(\mathfrak{p}+\mathfrak{q})+\eps.
\]
Thus
\[
b_i=\frac{\mathfrak{p}+\mathfrak{q}+3\eps}{N}.
\]
By Lemma~\ref{lem:rationalSeed}, the starting diagram
$\ov{s}^{3(i+1)}Q_{b_i}^0$ is well-defined, with triple
\[
(\bE_\la,\bE_1,\bE_0),
\]
where
\[
\bE_\la=(\mathfrak{p}-5\mathfrak{q},\mathfrak{p}-\mathfrak{q},2\eps),
\]
and
\[
\bE_0=(4\mathfrak{p}-23\mathfrak{q},\mathfrak{p}-4\mathfrak{q},\eps),
\qquad
\bE_1=(2\mathfrak{p}-11\mathfrak{q},\mathfrak{p}-2\mathfrak{q},-\eps).
\]

Since $\eps t_\la=2$, the $x$-mutation gives
\[
x^n\ov{s}^{3(i+1)}\Tt_0
=(\bE_\la,\bE_{n+1},\bE_n),
\qquad
\bE_{n+2}=2\bE_{n+1}-\bE_n.
\]
Consequently,
\[
\bE_n
=(p_n,q_n,t_n)=\bigl(
4\mathfrak{p}-23\mathfrak{q}+2n(6\mathfrak{q}-\mathfrak{p}),\,
\mathfrak{p}-4\mathfrak{q}+2n\mathfrak{q},\,
\eps(1-2n)
\bigr).
\]
Since $\mathfrak{p}/\mathfrak{q}>3+2\sqrt2>4$, we have $q_n>0$
for every $n\geq0$.

We first check that every $x$-mutation is well-defined.
Assume that $Q_n$ is well-defined; its triple is $(\bE_\la,\bE_{n+1},\bE_n)$. Since $6q_n-p_n=(2n+2)\mathfrak{p}-\mathfrak{q}>0$, Lemma~\ref{lem:barsCriterion} applies. Writing $(d,m)$ for the degree coordinates of $S^{-1}\bE_{n+1}$, we have $m-b_id=-\bigl((m_{n+1}-q_{n+1})-b_i(d_{n+1}-3q_{n+1})\bigr)$, and substituting the formulas for $\bE_{n+1}$ and $b_i$
into the degree formulas gives
\[
\eps\bigl((m_{n+1}-q_{n+1})
-b_i(d_{n+1}-3q_{n+1})\bigr)
=\frac{4(\mathfrak{p}-\mathfrak{q})}{N}>0.
\]
Hence $\eps(m-b_id)<0$, so by Lemma~\ref{lem:barsCriterion}(ii) the next $x$-mutation is well-defined. This proves
inductively that every $Q_n$ is well-defined.

We now compute the limiting triangle.
Write $Q_n=OX_nV_nY_n$, with
\[
X_n=(A_n,0),\qquad Y_n=(0,B_n).
\]
By \eqref{eq:algGeom}, $Q_n=Q_{b_i}(x^n\ov{s}^{3(i+1)}\Tt_0)$, and substituting the triple coordinates into Definition~\ref{def:Qb} and \cite[Lemma~2.3.3]{M2} gives
\[
A_n=\frac{8\mathfrak{p}}{N}-\frac{4}{Nq_n},
\qquad
B_n=\frac{8\mathfrak{q}}{N},
\qquad
|V_nY_n|=\frac{4}{Nq_{n+1}},
\]
where $|V_nY_n|$ denotes affine length.
Since $q_n\to\infty$, we obtain
\[
A_n\longrightarrow A:=\frac{8\mathfrak{p}}{N},
\qquad
B_n=B:=\frac{8\mathfrak{q}}{N}.
\]
The direction of $V_nY_n$ is determined by the fixed
left entry $\bE_\la$. Hence its affine length tending
to zero implies that $V_n-Y_n\to0$.
Therefore $Q_n$ converges to the triangle with
vertices $O$, $(A,0)$, and $(0,B)$.

The recursion for $(\mathfrak{p}_i,\mathfrak{q}_i)$ gives
\[
\mathfrak{p}^2-6\mathfrak{p}\mathfrak{q}+\mathfrak{q}^2=1,
\]
so $(\mathfrak{p}+\mathfrak{q})^2-1=8\mathfrak{p}\mathfrak{q}$. Consequently,
\[
1-b_i^2
=\frac{N^2-(\mathfrak{p}+\mathfrak{q}+3\eps)^2}{N^2}
=\frac{64\mathfrak{p}\mathfrak{q}}{N^2}
=AB.
\]
Moreover,
\[
\frac AB=\frac{\mathfrak{p}}{\mathfrak{q}}=\acc(b_i).
\]
Thus the limiting triangle gives a full filling
at $\acc(b_i)$.

Finally, $q_{n+1}>q_n$, so the formula for $A_n$
shows that $A_{n+1}>A_n$ for every $n$.
Hence the standardized polygons $Q_n$ are
pairwise distinct.
\end{proof}

 \bibliographystyle{plain}
\bibliography{references}

@article{MPW,
  title={A classification of infinite staircases for {H}irzebruch surfaces},
  author={Magill, Nicki and Pires, Ana Rita and Weiler, Morgan},
  journal={Journal of Topology},
  volume={18},
  number={1},
  pages={e70017},
  year={2025},
  publisher={Wiley Online Library}
}

@article{REU,
  title={Four-periodic infinite staircases for four-dimensional polydisks},
  author={Farley, Caden and Holm, Tara S and Magill, Nicki and Schroder, Jemma and Wang, Zichen and Weiler, Morgan and Zabelina, Elizaveta},
  journal={Involve, a Journal of Mathematics},
  volume={18},
  number={1},
  pages={25--78},
  year={2025},
  publisher={Mathematical Sciences Publishers}
}

@article{MMW,
  title={Staircase patterns in {H}irzebruch surfaces.},
  author={Magill, Nicki and McDuff, Dusa and Weiler, Morgan},
  journal={Commentarii Mathematici Helvetici},
  volume={99},
  number={3},
  pages={437--555},
  year={2024}
}

@article{MM,
  title={Staircase symmetries in {H}irzebruch surfaces},
  author={Magill, Nicki and McDuff, Dusa},
  journal={Algebraic \& Geometric Topology},
  volume={23},
  number={9},
  pages={4235--4307},
  year={2023},
  publisher={Mathematical Sciences Publishers}
}

@article{M1,
  author  = {Magill, Nicki},
  title   = {Unobstructed embeddings in Hirzebruch surfaces},
  journal = {Journal of Symplectic Geometry},
  volume  = {22},
  number  = {1},
  pages   = {109--152},
  year    = {2024},
  doi     = {10.4310/JSG.2024.v22.n1.a3}
}

@article{McDuffSiegel2025,
  author  = {McDuff, Dusa and Siegel, Kyler},
  title   = {Singular algebraic curves and infinite symplectic staircases},
  journal = {Inventiones Mathematicae},
  volume  = {242},
  pages   = {387--459},
  year    = {2025},
  doi     = {10.1007/s00222-025-01359-4}
}

@article{FibStairs,
  author  = {McDuff, Dusa and Schlenk, Felix},
  title   = {The embedding capacity of {$4$}-dimensional symplectic ellipsoids},
  journal = {Annals of Mathematics},
  volume  = {175},
  number  = {3},
  pages   = {1191--1282},
  year    = {2012},
  doi     = {10.4007/annals.2012.175.3.5}
}

@misc{CGMM,
  author       = {Cristofaro-Gardiner, Dan and Magill, Nicki and McDuff, Dusa},
  title        = {Curvy points, the perimeter, and the complexity of convex toric domains},
  howpublished = {arXiv:2506.23498},
  year         = {2025}
}

@article{CGHMP,
  title={On infinite staircases in toric symplectic four-manifolds},
  author={Cristofaro-Gardiner, Dan and Holm, Tara S and Mandini, Alessia and Pires, Ana Rita},
  journal={Journal of Differential Geometry},
  volume={129},
  number={2},
  pages={335--413},
  year={2025},
  publisher={Lehigh University}
}

@article{CV,
 title = {Full ellipsoid embeddings and toric mutations}, abstractNote = {Abstract This article introduces a new method to construct volume-filling symplectic embeddings of 4-dimensional ellipsoids by employing polytope mutations in toric and almost toric varieties. The construction uniformly recovers the full sequences for the Fibonacci Staircase of McDuff–Schlenk, the Pell Staircase of Frenkel–Müller and the Cristofaro-Gardiner–Kleinman Staircase, and adds new infinite sequences of ellipsoid embeddings. In addition, we initiate the study of symplectic-tropical curves for almost toric fibrations and emphasize the connection to quiver combinatorics.}, journal = {Selecta Mathematica}, volume = {28}, number = {3}, author = {Casals, Roger and Vianna, Renato}, year={2022} }

@incollection{ICERM,
  author    = {Bertozzi, Maria and Holm, Tara S. and Maw, Emily and McDuff, Dusa and Mwakyoma, Grace and Pires, Ana Rita and Weiler, Morgan},
  title     = {Infinite staircases for {H}irzebruch surfaces},
  booktitle = {Research Directions in Symplectic and Contact Geometry and Topology},
  editor    = {Acu, Bahar and Cannizzo, Catherine and McDuff, Dusa and Myer, Ziva and Pan, Yu and Traynor, Lisa},
  series    = {Association for Women in Mathematics Series},
  volume    = {27},
  pages     = {47--157},
  publisher = {Springer, Cham},
  year      = {2021},
  doi       = {10.1007/978-3-030-80979-9_2}
}

@article{usher,
  title={Infinite staircases in the symplectic embedding problem for four-dimensional ellipsoids into polydisks},
  author={Usher, Michael},
  journal={Algebraic \& Geometric Topology},
  volume={19},
  number={4},
  pages={1935--2022},
  year={2019},
  publisher={Mathematical Sciences Publishers}
}

@misc{M2,
  author       = {Magill, Nicki},
  title        = {Quadrilateral mutations and symplectic embeddings},
  howpublished = {arXiv:2606.12729},
  year         = {2026}
}

@incollection{S,
	author = {Margaret Symington},
	booktitle = {Topology and Geometry of Manifolds},
	editor = {Gordana Mati{\'c} and Clint McCrory},
	pages = {153-208},
	publisher = {American Mathematical Society},
	series = {Proceedings of Symposia in Pure Mathematics},
	title = {Four dimensions from two in symplectic topology },
	volume = {71},
	year = {2003}}

@book{Evans, place={Cambridge}, series={London Mathematical Society Student Texts}, title={Lectures on Lagrangian Torus Fibrations}, publisher={Cambridge University Press}, author={Evans, Jonny}, year={2023}, collection={London Mathematical Society Student Texts}}

@article{LS,
author = {Naichung Conan Leung and Margaret Symington},
title = {{Almost toric symplectic four-manifolds}},
volume = {8},
journal = {Journal of Symplectic Geometry},
number = {2},
publisher = {International Press of Boston},
pages = {143 -- 187},
year = {2010},
}

@misc{MarkovStairs,
  author       = {Adaloglou, Nikolas and Brendel, Jo{\'e} and Evans, Jonny and Hauber, Johannes and Schlenk, Felix},
  title        = {Markov staircases},
  howpublished = {arXiv:2509.03224},
  year         = {2025}
}
\end{document}